\documentclass[11pt,a4paper,english]{amsart}
\usepackage{amssymb,amsmath,amsthm,babel,bm}
\usepackage{graphicx,setspace}
\newtheorem{theorem}{Theorem}[section]
\newtheorem{proposition}[theorem]{Proposition}
\newtheorem{definition}[theorem]{Definition}
\newtheorem*{definition*}{Definition}
\newtheorem{lemma}[theorem]{Lemma}
\newtheorem*{lemma*}{Lemma}
\newtheorem{remark}[theorem]{Remark}

\newtheorem{corollary}[theorem]{Corollary} 
 
\newtheorem*{remark*}{Remark}

\newtheorem{claim}[theorem]{Claim}
\newtheorem{fact}[theorem]{Fact}
\newcommand{\Pt}{P}%A point

\newcommand{\N}{{\mathbb N}}
\newcommand{\Z}{{\mathbb Z}}
\newcommand{\Q}{{\mathbb Q}}
\newcommand{\C}{{\mathbb C}}

\newcommand{\Cu}{{\mathcal C}}
\newcommand{\F}{{\mathbb F}}

\newcommand{\Qb}{\overline{\mathbb Q}}

\newcommand{\A}{{\mathbb A}}
\renewcommand{\P}{{\mathbb P}}
\newcommand{\Gm}{{\mathbb G}_{\rm m}}
\newcommand{\codim}{{\rm codim}}
\renewcommand{\div}{{\rm div}}

\newcommand{\ord}{{\rm ord}}
\renewcommand{\d}{{\rm d}}
\newcommand{\ie}{{\it i.\,e.\ }}
\newcommand{\e}{{\varepsilon}}
\newcommand{\cf}{{\it cf}~}

\newcommand{\rhog}{{\boldsymbol\rho}}
\newcommand{\alphag}{{\boldsymbol\alpha}}
\newcommand{\betag}{{\boldsymbol\beta}}
\newcommand{\gammag}{{\boldsymbol\gamma}}
\newcommand{\detag}{{\boldsymbol\delta}}
\newcommand{\omegag}{{\boldsymbol\omega}}

\renewcommand{\v}{{\bf v }}
\newcommand{\wg}{{\bf w }}
\newcommand{\xg}{{\bf x }}
\newcommand{\gb}{{\bf g }}
\newcommand{\fb}{{\bf f }}

\newcommand\DIM{{\smallskip\noindent{\bf Proof.}\quad}}
\newcommand\CVD{\begin{flushright}$\square$\end{flushright}
\vskip 0.2cm\goodbreak}
\def\newatop#1#2{\genfrac{}{}{0pt}{}{#1}{#2}}

\numberwithin{equation}{section}

\title{Bounded Height in Pencils of Finitely Generated Subgroups.}
\author{F. Amoroso,  D. Masser and U. Zannier}
\begin{document}
\maketitle

\begin{abstract}  
In this paper we prove a general bounded height result for specializations in finitely generated subgroups varying in families which complements and sharpens the toric Mordell-Lang Theorem by replacing {\it finiteness} by {\it emptyness}, for the intersection of varieties and subgroups, all moving  in a pencil, except for bounded height values of the parameters (and excluding identical relations).

More precisely, an instance of the results is as follows. Consider the torus scheme ${{\Gm^r}_{/\Cu}}$ over a curve $\Cu$ defined over $\Qb$, and let $\Gamma$ be a subgroup-scheme generated by finitely many sections (satisfying some necessary conditions). Further, let $V$ be any subscheme. Then there is a bound for the height of the points $P\in\Cu(\Qb)$  such that,   for some  $\gamma\in\Gamma$ which does not generically lie in $V$,   $\gamma(P)$ lies in the fiber $V_P$.

We further  offer some direct  diophantine applications, to illustrate once again that the results implicitly contain information   absent from the previous bounds in this context.
%\noindent{\bf Mathematics Subject Classification:} %11G5O (Primary), 11Jxx (Secondary).
\end{abstract}  

\section{Introduction}

Let $\Cu$ be a projective smooth curve defined over $\Qb$, with function field  denoted $\F:=\Qb(\Cu)$. 
In 1999 Bombieri and the second and third authors \cite{Bo-Ma-Za} proved a bounded height result 
in the multiplicative group $\Gm^r$. 
\begin{theorem}[\cite{Bo-Ma-Za}, Theorem 1'] 
\label{BoMaZa}
Let $\Gamma\subset\Gm(\F)$ be a finitely generated subgroup of non-zero rational functions on $\Cu$ such that the only constants in $\Gamma$ are roots of unity. Then the height of the points $\Pt\in\Cu(\Qb)$, such that  for some $\xg\in\Gamma\setminus\{1\}$  we have $\xg(P)=1$, is bounded above.
\end{theorem}
A significant special case is 
$$
t^n(1-t)^m=1
$$
provided only $n$, $m$ are not both zero.\\

In this paper we prove a general bounded height result for specializations in finitely generated subgroups varying in families. This vastly  extends the previously treated constant case and involves entirely different, and  more delicate, techniques. 

Before stating our results, note that Theorem \ref{BoMaZa} may be also phrased as a kind of toric analogue of Silverman's Specialization Theorem (\cite{Si}, Theorem C). To illustrate this link, let us  consider the `trivial' family ${{\Gm}_{/\Cu}}:=\Gm\times\Cu$ and  the sections $\gamma_i:\Cu\to\Gm\times \Cu$, given by $P\mapsto (g_i(P),P)$, where $g_i$ are generators for $\Gamma$, independent modulo constants. The above conclusion then means that {\it the set of points $P$ where the values of the sections are multiplicatively dependent has bounded height}.

Now, multiplicative dependence at $P$ means that some nontrivial monomial attains the value $1$ at $P$. Then, rather more generally, given a constant family $\pi: {{\Gm^r}_{/\Cu}}:=\Gm^r\times\Cu\to \Cu$, a subvariety $V\subset \Gm^r\times\Cu$ and sections $\gamma_i:\Cu\to {{\Gm^r}_{/\Cu}}$ generating a group $\Gamma$, we may ask  the following\\

\noindent{\bf Question}: {\it What can be said about those points $P$ such that some nontrivial  element of the group $\Gamma$  when specialized at $P$ lies on the fiber $V_P=\pi^{-1}(P)\subset \Gm^r$}.\\

The previous situation  is obtained in the very special case when $V_P$ is constantly equal to  the origin, \ie $V=$ origin$\times\Cu$. By contrast, we stress  that  here neither $V$ nor the $V_P$   are assumed to have any kind of group structure; especially this feature heavily prevents the previously known proof-pattern to apply.\\

The present paper offers in a sense a complete solution to this issue, proving that on the appropriate assumptions we have generally bounded height for {\it any} proper family of subvarieties. In particular, the said intersection is empty except for a `sparse' set of points.
  
For simplicity, we  phrase this conclusion in the language of Theorem \ref{BoMaZa}. Namely, we consider a power $\Gm^r$ of the multiplicative algebraic group and we let $V$ be a subvariety of $\Gm^r$ defined over $\F$;  so we may view $V$ as a family of varieties parameterised by $\Cu$. Then we denote by $V_P$, for almost all $P\in\Cu$, a specialized variety defined e.g. by specializing at $P$ a given system of defining equations for $V$. Given a subgroup $\Gamma$ of $\Gm^r$ defined over $\F$ we say that $\Gamma$ is constant-free if its image $\Gamma'$ by any surjective homomorphism $\Gm^r\rightarrow\Gm$ satisfies the assumption  $\Gamma'\cap\Qb^*=\Gamma'_{\rm tors}$ of Theorem~\ref{BoMaZa}.

With such notation, we have the following uniform complement to the toric case of Lang's conjecture. 

\begin{theorem}
\label{specialization}
Let $\Gamma\subset\Gm^r(\F)$ be a finitely generated constant-free subgroup and let $V$ be a subvariety of $\Gm^r$  defined over $\F$. Then the height of the points $\Pt\in\Cu(\Qb)$, such that for some $\xg\in\Gamma\setminus V$ the value $\xg(P)$ is defined and lies in $V_P$, is bounded above.
\end{theorem}

Consider a generic situation when $\Gamma\cap V$ is empty and $V_P$ does not contain a coset of positive dimension. The Mordell-Lang Theorem tells us that $\Gamma_P\cap V_P$ is finite for all $P$. Theorem~\ref{specialization} gives the following complement: for $P$ of large height $\Gamma_P\cap V_P$ is {\it empty}.\\

As a consequence, we recover Theorem~\ref{BoMaZa} of \cite{Bo-Ma-Za}, taking $V=\{1\}$. An other known example is with $\Cu$ the affine line, $\F=\Qb(t)$ and $\Gamma\subset\Gm^2$ the subgroup generated by $(t,1-t)$. Let $V$ be the hypersurface of $\Gm^2$ defined by the equation $\alpha_1x_1+\alpha_2x_2=1$; here we obtain bounded height for 
\begin{equation}
\label{beuk}
\alpha_1t^n+\alpha_2(1-t)^n=1
\end{equation}
unless $n =1$ and $\alpha_1=\alpha_2=1$, a result of Beukers \cite{Be}.

Theorem~\ref{specialization} has also several entirely new consequences. For instance, choose $\Cu$ the affine line as above, $\Gamma\subset\Gm^3$ generated by $(t,1,1)$, $(1-t,1,1)$, $(1,t,1)$ and $(1,1,1+t)$. With $V$ of equation $x_1+x_2+x_3=1$ we get bounded height for the solutions of 
\begin{equation}
\label{amzex}
t^n(1-t)^m+t^l+(1+t)^p=1,
\end{equation}
with no proviso on $n,m,l,p$. More generally, any equation
\begin{equation}
\label{bmzbeuk}
\alpha_1M_1+\cdots+\alpha_sM_s=1
\end{equation}
for algebraic $\alpha_1,\ldots,\alpha_s$ and monomials $M_1\ldots,M_s$ in fixed algebraic functions of $t$ usually implies that the absolute height of $t$ is bounded above independently of $M_1,\ldots,M_s$. We can also think of allowing $\alpha_1,\ldots,\alpha_s$ in \eqref{bmzbeuk} to be not necessarily fixed, for example as fixed algebraic functions of $t$, or even as numbers varying at most subexponentially in the exponents occurring in $M_1,\ldots,M_s$.\\

\noindent{\bf Addendum to Theorem~\ref{specialization}.}\quad For a subvariety $V$ defined over the constant field $\Qb$ (\ie we have a  `trivial' family with $V_P=V$ for all $P$), the conclusion of Theorem~\ref{specialization} still holds for a subgroup $\Gamma$ which is not necessarily constant-free, but such that $\Gamma/\Gamma\cap\Gm^r(\Qb)$ is of rank $1$.\\

Note however that some assumption on $\Gamma$ is needed. As a non-trivial example, we may take as $\Cu$ the affine line as above, and $\Gamma\subset\Gm^2$ the subgroup generated by $\gamma_1=(t,1)$ and $\gamma_2=(1,2t)$ (note that these vectors are multiplicatively independent modulo constants). Let $V$ be the hypersurface of $\Gm^2$ defined by the equation $x_1+x_2=0$. Then for $n\in\N$ the element $\xg^{(n)}=\xg^{(n)}(t)=\gamma_1^{n+1}\gamma_2^n=(t^{n+1},2^nt^n)\in\Gamma$ is not generically in $V$ but its specialization at $t=-2^n$ is. Of course the image of $\Gamma$ under the isogeny sending $(x,y)$ to $xy$ contains 2.\\

Choosing $V$ a subvariety defined over the constant field $\Qb$, we obtain a bounded height result for certain `unlikely intersections'. 

\begin{corollary}
\label{unlikely}
Let $\Cu\subseteq\Gm^r$ be a curve and let $V\subseteq\Gm^r$ be a subvariety, both defined over $\Qb$. Then the height of the points $\Pt\in\Cu(\Qb)$, such that there exists an integer $n$ with\footnote{We denote as usual by $[n]\colon\Gm^r\rightarrow\Gm^r$ the morphism $\xg\mapsto\xg^n$ of  multiplication  by $n\in\Z$.} $[n]\Cu\not\subseteq V$ and $[n]\Pt\in V(\Qb)$, is bounded above.
\end{corollary}

Indeed, take $\Gamma=\langle(g_1,\ldots,g_r)\rangle$ where $g_1,\ldots,g_r$ are coordinate functions on $\Cu$, and let $n$ be an integer such that $[n]\Cu\not\subseteq V$. Then $\gamma:=(g_1^n,\ldots,g_r^n)\in\Gamma\setminus V$. By the addendum of Theorem~\ref{specialization}, the height of the points $\Pt\in\Cu(\Qb)$ such that $[n]\Pt\in V(\Qb)$ is bounded independently of $n$.\\

Note that by the Skolem-Mahler-Lech theorem, for a fixed point $P_0\in\Gm^r$ the set of integers $n$ with $[n]P_0\in V$ is a union of a finite number of points and arithmetic progressions, and is `usually' finite (if the set is not finite, $V$ must contain the Zariski-closure of a set $[nq+s]P_0$, $n\in\N$, for  suitable integers $q\neq 0, s$, and in particular must contain a coset of an algebraic subgroup of positive dimension unless $P_0$ is torsion). As in a comment above, the corollary says that  when we move $P_0$ along a curve $\Cu$ the corresponding set is `usually' empty, except for points $P_0\in\Cu(\Qb)$ of bounded height.

For more remarks on Theorem \ref{specialization} and Corollary~\ref{unlikely}, see section~\ref{further}.\\

Let us now discuss proofs. Those of \cite{Bo-Ma-Za} may be said to rely on some form of simultaneous diophantine approximation to rational numbers. But the proof of Beukers for (\ref{beuk}) is very different. The starting point is an explicit construction of Pad\'e approximants to powers of the linear polynomial $1-t$, which goes back to \cite{Be-Ti} and involves hypergeometric polynomials. This leads to explicit identities of the form $A(t)t^n+B(t)(1-t)^n=C(t)$ (suitable for (\ref{beuk}) above). 

In our general situation, or even just for $\gamma_1^n+\gamma_2^n=1$, we cannot hope to use such an explicit construction; just to mention one indication in this respect, results of Bombieri-Cohen (see \cite{Bo-Co}) suggest that in general  the coefficients of these Pad\'e approximants have a height growing much faster (\ie like $\exp (cn^2)$) than in the case of Beukers (when it grows like $\exp (cn)$), and this would destroy the basic  estimates necessary for the method to go through.

We use instead Thue's Method for avoiding such explicit constructions; this involves  divided derivatives, Siegel's Lemma and a zero estimate based on Wronskians. However,  to deal with certain unexpected vanishings we have to introduce a quite intricate descent, whose structure  is different with respect to  other investigations using Wronskians. This allows us to prove the following explicit version of a special case of Theorem~\ref{specialization}.

\begin{theorem}
\label{main}
Let $r\geq 2$ and $f_1,\ldots,f_r\in\F$ be non-zero rational functions such that $f_i/f_j$ is non-constant for some $i$ and $j$. Then there exists a positive real number $C$ depending only on $f_1,\ldots,f_r$, having the following properties. Let $\alphag=(\alpha_1:\cdots:\alpha_r)\in\P^{r-1}(\Qb)$.  
Consider, for a natural number $n$, a solution $P\in \Cu(\Qb)$ of the equation
\begin{equation}
\label{main-equation}
\alpha_1 f_1(\Pt)^n+\cdots + \alpha_r f_r(\Pt)^n =0.
\end{equation}
Then, if $n\geq C$ and if there are no proper vanishing subsums, we have 
\begin{equation}
\label{hP0}
h(\Pt)\leq \frac{r h(\alphag)}{n}+C.
\end{equation}
\end{theorem}

Our method provides even more explicit bounds: see Theorem~\ref{main-explicit} in Section~\ref{proofmain}, which is our main technical result. As also remarked by a referee we note that, given a fixed $\e> 0$, Theorem~\ref{main-explicit} allows to replace $r$ in formula~\eqref{hP0} by $r-1+ \e$ for fixed arbitrarily small $\e > 0$, by allowing $C$ to depend mildly on $\e$. Then remark~\ref{comments}, iii) shows that now this is sharp.

We finally remark that the assumption that $n$ is sufficiently large is necessary: we may have $\alpha_1f_1^n+\cdots + \alpha_r f_r^n=0$ identically  for some small values of $n$. Reciprocally, our results implies that $\alpha_1f_1^n+\cdots + \alpha_r f_r^n$ cannot in fact be zero for {\it infinitely many} $n$. This last fact  may be also derived directly on using  the main result of~\cite{Br-Ma}  (\ie  a general form of the $abc$-inequality over function fields).\footnote{It is no coincidence that the proof of this uses Wronskians, which also appear in our arguments; on the other hand, no consideration of heights appears in~\cite{Br-Ma}.}  \\

Of course the use of Thue's Method in classical diophantine approximation is well-known to lead to results which are usually not effective. By contrast all the results of this paper are effective. For example with $t^n+(1-t)^n+(1+t)^n=1$ (and $n \geq 0$) Denz \cite{Denz} gets $H(t) \leq 2^{856}$ for the standard (non-logarithmic, absolute) height.\\

\subsection{\tt Plan of the paper }
\label{plan}
In Section~\ref{auxi} we first define a simple notion of arithmetic height on rational functions and we prove some elementary estimates for this height, including some explicit bounds for the height of divided higher derivatives. Next we have two lemmas: Lemma~\ref{siegel-substitute} which is an auxiliary construction through Siegel's Lemma, and Lemma~\ref{Wronskiano} 
which constructs by specialization some approximation forms. This last lemma is a fundamental tool of the inductive proof in section~\ref{proofmain}.

In Section~\ref{strategy} we state Theorem~\ref{main-explicit}, which is a more precise and detailed version of Theorem~\ref{main}. Then, we first illustrate the strategy of the proof in a test-run case. For the general proof, we need to overcome a technical obstacle related to non-vanishing. At this stage we only state a Claim and we deduce the Theorem in full generality from it. 

Section~\ref{proofmain} is devoted to a proof of the Claim which relies on a somewhat intricate descent argument.

The final Section~\ref{coro} is devoted to the deduction of Theorem~\ref{specialization} from Theorem~\ref{main-explicit}, namely the case of arbitrary rank from the case of rank one. This is done using simultaneous diophantine approximation applied to the exponents of the group-generators.
\\ \\

\noindent{\bf Aknowledgements.} We thank the referees for their accurate consideration of the paper and several valuable suggestions. The first author thanks Bruno Angl\`es who draw his attention to Mirimanoff Polynomials.

\section{Applications and further remarks}
\label{app}

\subsection{\tt Two simple applications}\label{examples} We show by means of two  simple examples  that our results, actually already  very special cases of them, are capable of applications to  diophantine issues, recovering certain finiteness statements. The result achieved in the first example is known and may be obtained  by a variety of techniques, but here we reach  it directly as a consequence of  the above corollary. The second example appears to be new. Probably it can be sharpened, but our aim here is merely to illustrate possible applications, not to develop them in depth. For this same reason we shall be sketchy in the arguments.

\medskip

{\it Example 1: A family of Thue's equations}. Consider the Thue's equation   
$$
x^3-(t^3-1)y^3=1,
$$
where $t>1$ is an integer, to be solved in integers $x,y$.

Our results easily imply that {\it there are only finitely many integers $t$ such that the equation has an integer solution with $y\neq 0,1$.} 

\medskip

In fact, let  $u=\root 3\of {t^3-1}$, e.g. the real determination. It is very easy to see that the group of units of the ring $\Z[u]$ is  generated by $\pm 1$ and  the unit $t-u$. 

Let $(x,y)\in\Z^2$  be a solution, so $x-uy\in\Z[u]$ is a unit of norm $1$, whence   $=(t-u)^n$ for some integer $n$,  and  on taking conjugates over $\Q$ we have 
$$
(t-u)^n+\omega(t-\omega u)^n+\omega^2(t-\omega^2u)^n=0,
$$
where $\omega$ is a primitive  cube root of 1.

If  $n>1$, one may easily show that  the left side does not vanish identically. Hence, by the corollary (with $\Cu$ a curve with function field $\Qb(z,{\root 3\of {z^3-1}})$ and $r=3$),   the height of $t$ is bounded   and finiteness follows.

This particular result (even that there are no $t$) has been known for nearly a century, but our method allows substantial generalizations; for example it may be applied to the context of a conjecture of Thomas~\cite{Th} about families of Thue equations, and extended to non-Thue equations such as
$$x^3-(t^3-1)y^3+3(t^3-1)xy+(t^3-1)^2=1$$
(which cannot be obtained from a Thue equation by an inhomogeneous linear transformation), and even to more variables such as
$$x^4+(4t^4-1)y^4+(4t^4-1)^2z^4+2(4t^4-1)x^2z^2-4(4t^4-1)xy^2z = 1.$$
\medskip

{\it Example 2: Zeros of polynomial  recurrences}.  Consider a linear  recurrence sequence $(u_n(t))_{n\in\N}$ of polynomials in $t$, defined by prescribing   polynomial  initial data (not all zero) $u_0(t),\ldots ,u_{r-1}(t)\in\Q[t]$  and  imposing 
$$
u_{n+r}=c_1(t)u_{n+r-1}+\cdots +c_r(t)u_n,\qquad n=0,1,\ldots ,
$$
where $c_i(t)$ are polynomials  with  coefficients  say in $\Q$. We assume for simplicity that the characteristic polynomial $Z^r-c_1(t)Z^{r-1}-\cdots -c_r(t)$ has no multiple roots in an algebraic closure of $\Q(t)$. 

We assert that: {\it The set of algebraic numbers $\xi$ such that for some $n$ we have $u_n(t)\neq 0$ but $u_n(\xi)=0$ has bounded height. In particular, there are only finitely many such $\xi$ having bounded degree over $\Q$.}\\

Again, this follows rather immediately from the corollary and Northcott's theorem (see \cite{Bo-Gu}), after expressing $u_n(t)$ as a linear combination of $n$-th powers of the {\it roots} of the characteristic polynomial of the recurrence.

A rephrasing of the last conclusion is  that {\it  for every given $D$ there are only finitely many monic polynomials in $\Q[t]$ of degree $\le D$ which may divide some $u_n(t) \neq 0$}. 

In several cases we may further sharpen this,  to say something about the irreducible  factors. Assume  for instance that  $r\ge 3$ and that the characteristic polynomial of the recurrence has Galois group $S_r$ over $\bar\Q(t)$ (which is the `generic' case).\footnote{A doubly transitive Galois group would suffice for the sharpened conclusion. However it is beyond the scope of these examples to push the analysis further.} 
Under this assumption, we have:

\medskip

{\it  Apart from a finite set of polynomials, the degrees of the irreducible factors (over $\Q$) of the polynomials $u_n$ tends to infinity with $n$.}
 
 \medskip

To prove this claim note first that the assumption implies in particular that the recurrence is non-degenerate, meaning that no ratio of its roots is a root of unity. 
In turn,  by the Skolem-Mahler-Lech theorem,  this yields that  only finitely many $u_n$ may vanish identically. 

Suppose now that an algebraic number $\xi_0$ is a root of  $u_n$ for infinitely many $n$. Then, again by the Skolem-Mahler-Lech theorem  `several'  ratios of the specialised roots become roots of unity,  so the $n$-th powers of the roots of the recurrence collapse and vanishing of $u_n(\xi_0)$ gives rise to  linear relations among these roots of $1$.  If the order of these roots of unity is eventually unbounded, then  known theorems on torsion points on curves in $\Gm^r$ (see \cite{Za}) imply that the roots of the recurrence may be grouped in subsets of cardinality $\ge 2$ such that all the ratios of two roots in the same subset are multiplicatively dependent as algebraic functions; but the assumption on the Galois group  then easily yields a contradiction. On the other hand, if the order of the relevant roots of unity is bounded, then we find only finitely many $\xi_0$, proving the assertion. 

\medskip
 
 Note that some assumption as above is needed for this conclusion; an example in this direction is given by the Chebishev polynomials $T_n$ defined by $T_0(t)=2, T_1(t)=t$ and $T_{n+2}(t)=tT_{n+1}(t)-T_n(t)$ for $n\ge 0$. It turns out that, for odd $m$, $T_q$ divides $T_{mq}$ for all odd $q$, providing an example when the last conclusion is not true. Similarly for polynomials like $u_n(t):=T_n(t)T_{n+h}(t)$, where $r=4$ but the Galois group is too small. 

\subsection{\tt A relation with  Unlikely Intersections}\label{BH} The boundedness of the height in the set of $\Pt\in\Cu(\Qb)$ such that $[n]\Pt\in V(\Qb)$ is related also to  the context of Unlikely Intersections, and more precisely to  degenerate cases of the former Bounded Height Conjecture of Bombieri-Masser-Zannier~\cite{Bo-Ma-Za2}, nowadays a theorem of Habegger~\cite{Ha}, as in Example 1.3 of~\cite{Za}. 

Let $X\subseteq \Gm^r$ be a subvariety. Define $X^{\rm oa}$ as the complement in $X$ of the union of unlikely intersections of positive dimension, namely the components of some positive dimension $\delta>0$ of some intersection $X\cap B$, where $B$ is a translate of an algebraic subgroup and $\dim B\leq  \delta +\codim X-1$. By the former Habegger Theorem, the Weil height is bounded in the intersection of $X^{\rm oa}$  with the union of algebraic subgroups of dimension $\leq\codim X$. 

Let now $X=\Cu\times V\subset\Gm^r\times\Gm^r$. If $\Pt\in\Cu(\Qb)$ is such that $[n]\Pt\in V(\Qb)$, then $(\Pt,[n]\Pt)$ is in the intersection of $X$ with the algebraic subgroup $H_n=\{(\xg,\xg^n),\; \xg\in\Gm^r\}$ of dimension $\dim H_n=r\leq\codim X$. 
Thus it would be tempting to apply the former Bounded Height Conjecture = Habegger's Theorem. 
Nevertheless, such  result gives us no information here, since $X$ is {\it degenerate}: in the notation of \cite{Bo-Ma-Za2} we have $X^{\rm oa}=\emptyset$. 

This failure is not surprising,  because in the degenerate cases we cannot hope to have in general bounded height in the whole union $\bigcup_n (X\cap H_n)$. However, our Corollary  shows that,  at least, bounded height is recovered  in the projection of the whole union  to the first factor. \\

In the special case $V=\Cu$,  this kind of  problem  has been intensively studied by Bays and Habegger in~\cite{Ba-Ha}, who show (under suitable assumptions) the {\it finiteness} of the  set of $P\in\Cu$ such that  $[n]P\in\Cu$ for some $n\ge 2$, thus giving a partial answer to a question (with $[n]P\in\Cu'$ instead for another curve $\Cu'$) originated by A. Levin (see \cite{Za} notes of Chapter 1, p.39). As an important  tool  they first observe that the height in this set is bounded, except in the trivial case when $\Cu$ is a translate of a subtorus (op.cit., Lemma 6). This last result is a rather direct consequence  of a generalized Vojta's inequality due to R\'emond~\cite{Re2} (which however does not seem to apply when $\Cu' \neq \Cu$). This finiteness result provides further evidence of the usefulness of height bounds as in this paper. 

\subsection{\tt Abelian analogues?}\label{abelian} 
In this paper we limit ourselves to the toric case, but  analogous questions and statements can be naturally  formulated in the elliptic and abelian context, thus extending  Silverman's very setting.     

Note that Silverman's specialization results in  fact concern mainly abelian families, actually non-isotrivial, whereas algebraic  tori have only isotrivial\footnote{By `isotrivial' we mean that the family becomes trivial, \ie a product, after suitable extension of the base.} families; this introduces some differences  in the assumptions, as in results by Manin and Demjanenko, see \cite{Se}; for instance the assumption of mere {\it  independence} has to be  strengthened, as above,  in the  sense {\it modulo constants}, in order to obtain bounded height.

Further differences  with respect to the abelian case  are due to the fact that  heights in abelian varieties behave like quadratic forms, so  somewhat more `regularly' than in the toric case, when this lack may introduce difficulties in some parts of the proofs for the toric case.

A specific example  of what would be  the analogue of our main issue for  that case appears in the paper \cite{Ma-Za}, where we find  the pencil $\mathcal J$  of Jacobians $J_t$ of the curves     $H_t: y^2=x^6+x+t$, of generic genus $2$, parametrized by $t\in\A^1$; we also find  the section  $\sigma:\A^1\to\mathcal J$, obtained by defining   $\sigma(t)\in J_t$ as the class of the divisor  $\infty_+-\infty_-$,  difference of the two poles of $x$  on $H_t$ (let us forget here about the values for which $H_t$ has genus $<2$).   One can now consider the subvariety $V$ of $\mathcal J$ obtained as the union of $H_t$, embedding $H_t$ in $J_t$ e.g. through $\infty_+$.  

The question now is: {\it What can be said about the points $t_0\in\A^1$ such that $[n]\sigma(t_0)\in H_{t_0}$, for some $n=n(t_0)>3$?} (We exclude here $n=3$ because $[3]\sigma(t)\in H_t$ identically.)    Now even to prove that this set of $t_0$ (necessarily algebraic) has infinite complement in  $\Qb$ is far from being evident; the Appendix by Flynn to \cite{Ma-Za} in particular achieves this, and much more, providing nontrivial congruence conditions on the suitable $t_0$. However Flynn's method is not guaranteed to work generally, and moreover  the question of whether  the set of  these numbers has or not bounded height remains open. In fact,  we do not know  if and how the present methods can be adapted to the abelian context. It seems to us a rather  interesting issue to obtain even such a  hyperelliptic analogue. Or even a ``constant elliptic'' analogue; for example the points $P$ in $\Cu$ such that $[n]P \in \Cu'$ for some $n>1$, where $\Cu$, $\Cu'$ are defined by $x_1+x_2=1$, $x_1+x_2=c$ respectively in the product of $y_1^2=x_1^3+x_1+1$ and $y_2^2=x_2^3+x_2+1$ (we thank a referee for pointing out that when $\Cu'=\Cu$ this abelian version of Levin's question is covered by R\'emond's abelian analogue \cite{Re1} of  \cite{Re2}).

\subsection{\tt Further remarks }\label{further}
Let $\Cu\subseteq\Gm^r$ be a curve and let $V\subseteq\Gm^r$ be a subvariety, both defined over $\Qb$. Concerning the assumption $[n]\Cu\not\subseteq V$ in Corollary~\ref{unlikely}, a relevant issue is to detect the set of integers $n$ such that $[n]\Cu\subseteq V$. Now, this amounts to $[n]\xg\in V$ where the coordinates of $\xg$ are the restrictions to $\Cu$ of the coordinate functions on $\Gm^r$.  Classical results describe the Zariski closure of these $[n]\xg$ as a finite union of cosets of algebraic subgroups. Thus, if  $\Cu$ is not contained in any translate of a proper subtorus, no proper algebraic subgroup can contain a multiple $[m]\xg$ for $m\neq 0$, proving that  the said set is finite unless $V$ is the whole space\footnote{See also the paper \cite{SV} by Silverman and Voloch for more general finiteness results in this direction.}. Actually,  using for instance results in \cite{Br-Ma},   it is not difficult to reach directly this finiteness conclusion, moreover determining effectively the set.

A similar remark holds for Theorem \ref{specialization}: we may compute `effectively' the intersection $\xg\in\Gamma\cap V$, which however may be infinite in general. 

\section{Notation and Auxiliary Results.}
\label{auxi}

\subsection{Rational functions}

Given $f_1,\ldots,f_r\in\F$ not all zero, we put\footnote{Since we have chosen a smooth projective model of $\Cu$, the closed points over $\Qb$ of the curve correspond to the places of its function field $\F=\Qb(\Cu)$.} 
$$
\div(f_1,\ldots,f_r):=\sum_P \min_j\ord_P(f_j) P.
$$ 
Note that $\deg(\div(f_1,\ldots,f_r))=\deg(\div(gf_1,\ldots,gf_r))$ for any nonzero $g\in\F$.

Let us denote $d:=-\deg(\div(f_1,\ldots,f_r))$. We remark that $d\geq0$, since for example if $f_1 \neq 0$ we have $\sum_P \min_j\ord_P(f_j)\leq \sum_P\ord_P(f_1)=0$. Moreover $d>0$ if and only if some $f_i/f_j$ is non-constant. 

As a special case, let $f$ be a non-zero rational function on $\Cu$. We define as usual  its degree $d(f)$ as the the degree of the polar divisor $\div(f)_{\infty}=-\div(1,f)$. This is the geometric height of $(1:f)\in\P^1(\F)$. 

\medskip

{\tt An arithmetic height on $\F$}.   We define an arithmetic height $h(\cdot)$ of a rational function $f$ on $\Cu$ as follows. We choose once and for all a non-constant $t\in\F\backslash\Qb$. Let $F(X,Y)\in\Qb[X,Y]$ be the irreducible polynomial such that $F(t,f)=0$ (note that $F$ has degree at most $d(f)$ in $X$ and at most $d(t)$ in $Y$). 

\begin{definition} 
\label{def.height}  
For a function $f\in\F$, we define the height $h(f)$ as the projective Weil height of the vector of the coefficients of $F$.
\end{definition}

Clearly $h(1/f)=h(f)$. Also, this coincides with the affine height on $\Qb[t]$ (if $\Cu$ is the affine line and $f=P(t)$ is a polynomial then $F(X,Y)=P(X)-Y$.)\\

We shall need  the following elementary estimates for this height; we could not find suitable references in the standard literature, and indeed the very definition of $h(f)$ is not quite so standard.
\begin{lemma}
\label{height-stuff}
There is a constant $c$ depending only on $\Cu$ and $t$ with the following properties. Let $f\in\F$. Then we have 
\begin{itemize}
\item[(i)] $h(f^n) \leq nh(f)+cnd(f)$ for any positive integer $n$, 

\item[(ii)] $h(f') \leq c(h(f)+d(f)), d(f') \leq cd(f)$ for $f'=\d f/\d t$, 

\item[(iii)] $h({\rm tr}f) \leq h(f)+\log d(t), d({\rm tr}f)\leq d(t)d(f)$ for the trace  $tr$  from $\F$ to $\Qb(t)$, 

\item[(iv)] for $g\in\F$ we have
$$\max\{h(f+g),h(fg)\} \leq c(h(f)+h(g)+d(f)+d(g)),$$

\item[(v)] for any non-constant $s\in\F\backslash\Qb$ there is $C$ depending on $s$ (and $\Cu,t$) such that the height of $f$ with respect to $s$ is at most $C(h(f)+d(f))$. 

\end{itemize}
\end{lemma}
\DIM
We use resultants. Let $F(X,Y)$ be the irreducible polynomial for $f$ over $\Qb(t)$: $F(t,f)=0$. For (i) we can assume $d(f) \geq 1$, and we take the resultant of $F(X,Y)$ and $Y^n-Z$ with respect to $Y$ to get a non-zero polynomial $F_n(X,Z)$ with $F_n(t,f^n)=0$. Its degrees in $X,Z$ are at most $nd(f),d(t)$ respectively. And its height is at most $n(h(f)+c)$; here one must be careful to avoid a factorial in the number $N$ of terms in the Sylvester determinant, but it is easy to see that $N \leq (d(t)+1)^n2^{d(t)}\leq c^n$. Now $F_n$ might not be irreducible, but by well-known estimates the height of any factor is at most 
$$n(h(f)+c)+c(nd(f)+d(t)) \leq nh(f)+cnd(f);$$ 
and (i) follows.

A similar argument works with (v), now taking the resultant of $F(X,Y)$ and $S(X,Z)$ with respect to $X$, where $S(t,s)=0$. Also with (iv): now say $G(t,g)=0$ and then for $f+g$ we take the resultant of $F(X,Y)$ and $G(X,Z-Y)$ with respect to $Y$; then we do $g/f$ with $F(X,Y)$ and $G(X,ZY)$ and deduce $fg$ using $h(1/f)=h(f)$.

And (iii) is rather easy: if 
$$F(X,Y)=F_0(X)Y^e+F_1(X)Y^{e-1}+\cdots+F_e(X)$$ 
then ${\rm tr}f=-{d(t)\over e}{F_1(t) \over F_0(t)}$. Finally for (ii) we can also assume $d(f) \geq 1$ and then we note that
$$f'=-{F_0'(t)f^e+\cdots+F_e'(t) \over eF_0(t)f^{e-1}+\cdots+ F_{e-1}(t)}$$
so that we can use (iv).

\CVD

As at the beginning, we choose once and for all a system of Weil's functions associated to a divisor of degree $1$ and a corresponding height $h$ on $\Cu(\Qb)$. By a well-known result of N\'eron (see~\cite{Ne}), the height of $\Pt\in\Cu(\Qb)$ differs from $h(t(\Pt))/d(t)$ by an error term bounded by  a constant multiple of $1+h(\Pt)^{1/2}$.\\

We need the following functorial bound for the arithmetic height associated to values of $(f_1,\ldots,f_r)$, which is an easy consequence of ``Weil's Height Machine":

\begin{lemma}
\label{HeightMachine}
For $r \geq 2$ let $f_1,\ldots,f_r\in\F$ and $\Pt\in\Cu(\Qb)$, not a pole or a common zero of $f_1,\ldots,f_r$. Put $d:=-\deg\div(f_1,\ldots,f_r)$. Then the projective Weil height
$$
h(f_1(\Pt):\cdots:f_r(\Pt)) = d h(\Pt) + O(1+h(\Pt)^{1/2})
$$
where the implicit constant in the big-$O$ may depend on $f_1,\ldots,f_r$ but not on $\Pt$.
\end{lemma}
\DIM 
Let $E'=-\div(f_1,\ldots,f_r)$. We may assume $f_1,\ldots,f_r$ linearly independent over $\Qb$; indeed, we may select a basis, say $f_1,\ldots,f_s$ ($s\geq 2$), of the $\Qb$-vector space generated by $f_1,\ldots,f_r$ and observe  that $\div(f_1,\ldots,f_s)=\div(f_1,\ldots,f_r)$ and $h(f_1(\Pt):\cdots:f_r(\Pt)) = h(f_1(\Pt):\cdots:f_s(\Pt))+O(1)$. 

Given a divisor $D$, we denote by $h_D$ the height  on $\Cu(\Qb)$ associated to it, which is defined up  the addition of a bounded term: see~\cite{Hi-Si}, Part B for details. 

Let $\phi\colon\Cu\rightarrow\P^{r-1}$ be the morphism $\Pt\mapsto (f_1(\Pt):\ldots :f_r(\Pt))$ and $H$ be an hyperplane of $\P^{r-1}$. Then $\phi^*H$ is linearly equivalent to $E'$ (indeed, for $i=1,\ldots,r-1$, set $g_i=f_i/f_r$; then $\phi^*\{x_r=0\}=-\div(g_1,\ldots,g_{r-1},1)\sim E'$). 

Thus $h(f_1(\Pt):\cdots:f_r(\Pt))=h_{E'}(\Pt)+O(1)$. Here and in the rest of this proof, the big-$O$ depend on the divisors. 

We now apply Theorem B.5.9 of~\cite{Hi-Si}) (which goes back to N\'eron), taking as the ample divisor the divisor $D$ of degree $1$ such that  $h=h_D+O(1)$ and as the divisor equivalent to zero the divisor $E=dD-E'$. We obtain $h_E(\Pt)\leq c(1+h_D(\Pt)^{1/2})$. Moreover, by Theorem B.3.2 (d) of~\cite{Hi-Si}, $h_E=dh-h_{E'}+O(1)$. Thus 
\begin{align*}
h(f_1(\Pt):\cdots:f_r(\Pt))&=h_{E'}(\Pt)+O(1)\\
&=dh(\Pt)-dh_E(\Pt)+O(1)\\
&= dh(\Pt)+O(1+h(\Pt)^{1/2}).
\end{align*}
\CVD

 {\sc Uniformity.}  {\rm The question of how the implicit constants in the $O$-terms in the lemma depend on the functions is a subtle one, and  has been treated in number of papers, which control this dependence in data as the degree and heights of the functions. Here we can prove (at the expense of extra complication) the upper bound
$$
dh(P)+O\big((d+\max h(f_i))(1+h(P)^{1/2})\big)
$$
{with an implicit constant in the big-$O$ which depends only on $\Cu$ and $t$. If some sort of refined Height Machine could deliver the analogous lower bound, even for $r=2$, then it would imply at once some significant cases of our Theorem \ref{main}. For example with fixed different $F_1,\ldots,F_r$ in $\Qb[t]$ of degree $p \geq 1$ and fixed sufficiently general $\alpha_1,\ldots,\alpha_r$ in $\Qb$ we would have $d=np$ for
$$f(t)={\alpha_1F_1(t)^n+\cdots+\alpha_{r-1}F_{r-1}(t)^n \over \alpha_rF_r(t)^n}.$$
Thus when $f(t)=-1$ we would deduce
$$0=h(f(t)) \geq dh(t)-O\big(n(1+h(t)^{1/2})\big)$$
so $h(t)=O(1)$.

However such a lower bound is false in general, as the example $f(t)=(t-2^d)t^{d-1}-1$ with $f(2^d)=-1$ shows; the lower bound would be
$$dh(2^d)-O\big(d(1+h(2^d)^{1/2})\big)=d^2\log 2-O(d^{3/2}),$$
a contradiction for sufficiently large $d$.
}

\medskip

We shall need also some good bounds for the values of higher derivatives $f^{(l)}=(\d/\d t)^lf$; or rather those of the divided derivatives $\delta_lf=f^{(l)}/l!$ essential to the success of Thue's Method. Iterating part (ii) of Lemma \ref{height-stuff} does not suffice. In fact we have to consider certain monomial expressions whose curious weighting will soon be justified.

\begin{lemma}
\label{heisenstein}
For any $f\in\F$ there is $c$, depending only on $f$ and $t$, with the following property. Suppose $f$ and $t$ are regular at some $\Pt\in\Cu(\Qb)$ with $dt(\Pt) \neq 0$. For any non-negative integer $L$ let $\detag$ be any vector with components 
$$\delta_0f(\Pt)^{a_0}\delta_1f(\Pt)^{a_1}\cdots\delta_lf(\Pt)^{a_l}$$
for non-negative exponents satisfying
$$a_0+a_1+\cdots+a_l\leq L,~~a_1+\cdots+la_l\leq L.$$ 
Then the affine Weil height $h(\detag)$ is at most $cL(h(\Pt)+1)$. 

\end{lemma}
\DIM
This falls into the circle of Eisenstein-related ideas. We can assume that $f$ is not constant. With $F(t,f)=0$ as above and $\alpha_j=\delta_jf(\Pt)$ the power series $y=\sum_{j=0}^\infty \alpha_jx^j$ satisfies $F_0(x,y)=0$ with $F_0(X,Y)=F(t(\Pt)+X,Y)$. We may therefore apply Theorem 1 (p.162) of Schmidt's paper \cite{Sc1}. He needs a number field $k$ over which $F_0$ is defined. As he notes, the $\alpha_j$ lie in an extension $K$ of $k$ of relative degree at most the degree $e \geq 1$ of $F_0$ in $Y$; thus $[K:k]\leq d(t)$. We find for each valuation $v$ on $k$ some $A_v \geq 1$, with $A_v=1$ for all but finitely many $v$, such that
$$|\alpha_j|_w \leq A_v^{m+j}~~~~(j=0,1,2,\ldots)$$
for any valuation $w$ on $K$ over $v$, where $m\geq 1$ is the degree of $F_0$ in $X$; thus $m \leq d(f)$. Thus
$$|\alpha_0^{a_0}\alpha_1^{a_1}\cdots\alpha_l^{a_l}|_w\leq A_v^{m(a_0+a_1+\cdots+a_l)+(a_1+\cdots+la_l)}\leq A_v^{2mL}.$$
It follows for the non-logarithmic height
$$H(\detag)^{[K:\Q]} ~\leq~ \prod_vA_v^{2mL[K:k]}.$$
The $v$ are split into two sets $S_{\infty1},S_2$ with
$$\prod_{v\in S_{\infty1}}A_v ~\leq~ ((m+1)(e+1)\sqrt{e})^{(2e+1)[k:\Q]}H(F_0)^{2e[k:\Q]}~\leq~ (2H(F_0))^{c[k:\Q]},$$
$$\prod_{v\in S_2}A_v ~\leq~ (16m)^{11e^3[k:\Q]}H(F_0)^{(2e^3+2e)[k:\Q]}~\leq~ (2H(F_0))^{c[k:\Q]},$$
where $H(F_0)$ is still projective (and absolute). Thanks to the crucial (but nowadays natural) linear dependence on $[k:\Q]$ in the exponents we deduce $H(\detag)\leq (2H(F_0))^{4cmL}$.
Finally $H(F_0) \leq cH(t(\Pt))^m$, and the result we want follows by relating $h(t(\Pt))$ to $h(\Pt)$ as described above.

We remark that the extra precision of \cite{Sc1} (especially concerning the set $S_2$) is not really necessary for us; thus by putting harmless additional restrictions on the point $\Pt$ we could have got ourselves into the ``non-singular" situation, where the proofs are much easier (as for example in \cite{Bo-Gu} p.360).

\CVD

Given a divisor $D$ we denote by $L(D)$ the finite-dimensional $\Qb$-vector space 
$$
L(D)=\{f\in\F^*,\; \div(f)+D\geq 0\}\cup\{0\}.
$$ 
and by $l(D)$ its dimension. We shall need a good basis of $L(NQ)$ for fixed $Q$ and large $N$. It is convenient to talk also of $L(\infty Q)=\sum_{N=1}^\infty L(NQ)$ the vector space of $f\in\F$ which are regular outside $Q$.

\begin{lemma}
\label{good-basis}
For any $Q\in \Cu(\Qb)$ there is a positive integer $\Delta$ and real $c$ together with $g,g_0,g_1,\ldots,g_{\Delta-1}$ in $L(\infty Q)$, depending only on $Q$ and $t$, such that the following hold for any $N \geq 1$. 
\begin{itemize}
\item[(i)] We have $d(g)=\Delta$, and the elements
$$g_jg^k~~~(j=0,1,\ldots,\Delta-1,~~d(g_j)+kd(g) \leq N)$$
form a basis for $L(NQ)$, with
$$N-c \leq l(NQ) \leq N+1.$$

\item[(ii)] For any $f$ in $L(NQ)$ we have
$$f=\sum\alpha_{jk}g_jg^k$$
with affine height
$$h(\dots \alpha_{jk} \ldots) \leq c(h(f)+N),$$
as well as
$$h(f) \leq c(h(\dots \alpha_{jk} \ldots)+N).$$
\end{itemize}
\end{lemma}
\DIM
It is well known, for example by the Riemann-Roch Theorem, that as $f \neq 0$ varies over $L(\infty Q)$ the $-{\rm ord}_Q(f)$ (which are none other than the degrees $d(f)$) take all sufficiently large values. Let $\Delta \geq 1$ be the smallest positive value, attained by some $g$ in $L(\infty Q)$. If $\Delta=1$ we are done, as a standard argument of killing poles shows (which will be repeated below).

So we may and shall assume $\Delta \geq 2$. For $j=1,\ldots,\Delta-1$ pick $g_j$ in $L(\infty Q)$ with $n_j=-{\rm ord}_Q(g_j) \equiv j$ modulo $\Delta$ and also as small as possible; here $n_j > 0$ is automatic and even $n_j > \Delta$. We define $g_0=1$ and $n_0=0$. We show by induction on $N$ that these do the trick in (i).

Pick any $f \neq 0$ in $L(NQ)$, so that $n=-{\rm ord}_Q(f)\leq N$. If $n \equiv j$ modulo $\Delta$ $(j=0,1,\ldots,\Delta-1)$ then $n \geq n_j$ and we can find $\alpha$ in $\Qb$ with $f-\alpha g_jg^{(n-n_j)/\Delta}$ in $L((N-1)Q)$. As
$$d(g_j)+{n-n_j \over \Delta}d(g) = n \leq N$$
this shows by induction that the elements in (i) span $L(NQ)$. They are certainly linearly independent, as the
$$-{\rm ord}_Q(g_jg^k)=n_j+k\Delta ~~~\left(j=0,1,\ldots,\Delta-1,~~k=0,1,\ldots, \left[{N-n_j \over \Delta}\right]\right)$$
are all different (even for all $k$). In particular 
$$l(NQ)=\sum_{j=0}^{\Delta-1}\left(1+\left[{N-n_j \over \Delta}\right]\right),$$
which leads easily to the required estimates using $n_j \geq \Delta$ for $j>0$. So (i) is proved. Incidentally it is not difficult to estimate the constants so far solely in terms of the genus of $\Cu$.

For (ii) we note that this last argument even shows that $g_0,g_1,\ldots,g_{\Delta-1}$ are linearly independent over $\Qb(g)$. As $\Delta=d(g)=[\F:\Qb(g)]$ it follows that they form a basis of $\F$ over $\Qb(g)$. Now we can write 
$$f=\sum_{j=0}^{\Delta-1}g_jG_j$$ 
with $G_j=\sum_k\alpha_{jk}g^k$ in $\Qb[g]$. In the standard way we multiply by $g_0,g_1,\ldots,g_{\Delta-1}$ and take the trace from $\F$ to $\Qb(g)$. The resulting equations can be solved for $G_0,G_1,\ldots,G_{\Delta-1}$. It follows easily from Lemma \ref{height-stuff} (iii) (with $g$ not $t$),(iv),(v) that
$$\sum_{j=0}^{\Delta-1}h(G_j) \leq c(h(f)+N).$$
But the affine height in (ii) is at most the analogous sum with heights taken with respect to $g$. So the first of the two required inequalities follows with another appeal to Lemma \ref{height-stuff} (v). The second is similar but easier.

\CVD

Before we go further we record the following identity for divided derivatives. Namely
\begin{equation}
\label{francesco}
\delta_l(f^n)=f^{n-l}\sum_{\bf a}C({\bf a})\delta_0(f)^{a_0}\delta_1(f)^{a_1}\cdots\delta_l(f)^{a_l}
\end{equation}
where the sum is taken over all ${\bf a}=(a_0,a_1,\ldots,a_l)$ with non-negative coordinates satisfying
$$|{\bf a}|=a_0+a_1+\cdots+a_l=l,~~a_1+\cdots+la_l=l$$
(see earlier) and the $C({\bf a})$ are non-negative integers. Some version for undivided derivatives is attributed to the Blessed Francesco Fa\`a di Bruno (who even has $\phi(f)$ instead of $f^n$), but in this divided form we get an immediate proof by formally writing $\tilde f=\sum_{m=0}^\infty\delta_m(f)T^m$. Note that, formally, $\tilde f(x)=f(x+T)$.  

Taking this into account   one may then  go ahead by picking out the coefficient of $T^l$ in $\tilde f^n$. We take $\delta_0(f)=f$ in $a_0'$ of the factors $\tilde f$, and then $\delta_1(f)$ in $a_1$ of the factors, and so on. Then $a_0'+a_1+\cdots+a_l=n$ and $a_1+\cdots+la_l=l$ making it clear that $a_0=a_0'-(n-l)\geq 0$. We need also good estimates for the $C({\bf a})$, but it is similarly clear that their sum is majorized by the coefficient of $T^l$ in $(1+T+T^2+\cdots)^n=(1-T)^{-n}$, which is 
$$(-1)^l{-n \choose l}={n+l-1 \choose l} \leq 2^{n+l}$$
and in particular factorial-free.

Now comes our basic ``auxiliary polynomial".
\begin{lemma}
\label{siegel-substitute}
Let $f_1,\ldots,f_r$ be in $\overline{\Q}({\Cu})$, let $Q$ in $\Cu(\overline{\Q})$ be not a pole or zero of $f_1,\ldots,f_r$ also with $dt(Q) \neq 0$, and write as before $d=-\deg{\rm div}(f_1,\ldots,f_r)\geq 0$. Then there are $c_0,c$ depending only on $\Cu,f_1,\ldots,f_r,t$ and $Q$ with the following property. For any non-negative integers $n,M_1,\ldots,M_r$ with $S>M+nd+c_0$ for
$$S=M_1+\cdots+M_r,~~M=\max\{M_1,\ldots,M_r\}$$
define the ``Dirichlet exponent"
$$\varrho={M+dn \over S-M-dn-c_0}.$$ 
Then there are $A_1 \in L(M_1Q),\ldots,A_r \in L(M_rQ)$, not all zero and with heights at most $c(\varrho+1)(M+n)$, such that $A_1f_1^n+\cdots+A_rf_r^n=0$.
\end{lemma}
\DIM
By Lemma \ref{good-basis} (i) we can take
\begin{equation}\label{f}
A_i=\sum\alpha_{ijk}g_jg^k~~~~(i=1,\ldots,r)
\end{equation}
with $j,k$ satisfying
$$j=0,1,\ldots,\Delta-1,~~d(g_j)+kd(g) \leq M_i$$
and algebraic numbers $\alpha_{ijk}$ to be determined. We first find them such that if $\phi=\sum_{i=1}^rA_if_i^n \neq 0$ then
\begin{equation}\label{a}
{\rm ord}_Q\phi > T
\end{equation}
where the integer $T$ is nearly as large as linear algebra allows.

If $U$ is the number of unknowns $\alpha_{ijk}$ then their vector $\bm{\alpha}$ must lie in a certain subspace $V$ of $\Qb^U$. Here $U = \sum_{i=1}^rl(M_iQ)$ so Lemma \ref{good-basis}(i) gives
\begin{equation}\label{b}
S-c \leq U \leq S+r.
\end{equation}
We have $E=T+1$ equations, so the dimension $D$ of $V$ satisfies
\begin{equation}\label{c}
D \geq U-E \geq S-T-c_0.
\end{equation}
Thus we assume $T<S-c_0$ for solvability. But if $T$ is too near $S$ then as in Thue's Method we would lose control of the heights. To regain this we use the version of the Absolute Siegel Lemma proved by David and Philippon \cite{Da-Ph}; for example taking $\epsilon=1$ in the estimate at the bottom of page 523 we find non-zero $\bm\alpha\in V$ with 
\begin{equation}\label{d}
h(\bm\alpha) \leq {h(V) \over D}+{1 \over 2}\log D+1
\end{equation}
where $h(V)$ is the euclidean height. 

To estimate $h(V)$ we note that it is defined by certain equations, and by Hadamard's inequality for determinants (here the factorials don't matter) we get 
\begin{equation}\label{e}
h(V) \leq cE(\log E+h_{\rm eq})
\end{equation}
where
$h_{\rm eq}$ is an upper bound for the logarithmic euclidean height of each equation. These are $\delta_l\phi(Q)=0$ or more explicitly
$$
\sum_{i=1}^r\sum_{j,k}\alpha_{ijk}\beta_{ijkl}=0~~~~(l=0,1,\ldots,T)
$$
with $\beta_{ijkl}=\delta_l(g_jg^kf_i^n)(Q)$. Now 
$$\delta_l(g_jg^kf_i^n)=\sum\delta_s(g_j)\delta_p(g^k)\delta_q(f_i^n)$$
taken over all non-negative integers $s,p,q$ with $s+p+q=l$. By~\eqref{francesco}
$$\delta_p(g^k)=g^{k-p}\sum_{\bf a}C({\bf a})\delta_0(g)^{a_0}\delta_1(g)^{a_1}\cdots\delta_p(g)^{a_p},$$
$$\delta_q(f_i^n)=f_i^{n-q}\sum_{\bf b}C({\bf b})\delta_0(f_i)^{b_0}\delta_1(f_i)^{b_1}\cdots\delta_q(f_i)^{b_q}.$$
Now we see without difficulty thanks to Lemma \ref{heisenstein} that $h_{\rm eq} \leq c(T+M+n)$. Then (\ref{d}) and (\ref{e}) lead to
$$h(\bm\alpha) \leq c{(T+M+n)(T+1) \over S-T-c}+\log (S+c)$$
because $S-T-c_0 \leq D \leq U \leq S+r$ by (\ref{b}) and (\ref{c}).

We now choose $T$ so large that the condition (\ref{a}) forces after all $\phi=0$ in the sense that ${\rm ord}_Q\phi=\infty$. In fact (\ref{a}) holds also for $\tilde \phi=\phi/f_r^n$ because $f_r(Q) \neq 0$, and since $A_1,\ldots, A_r$ are in $L(MQ)$ it is easy to see that $d(\tilde \phi) \leq M+dn$. So $T=M+dn$ will do, leading to
$$h(\bm\alpha) \leq c\varrho(M+dn)+\log (S+c).$$
Finally Lemma \ref{good-basis} (ii) gets us to $h(A_i)$ by (\ref{f}); and then we use $S \leq rM$.
\CVD

\subsection{Orthogonal spaces and key lemma}
\label{orto}
~\\
Let $\wg$ be a fixed vector of $\C^n$ with all entries non-zero. Given a subset $\Lambda$ of $\{1,\ldots,r\}$ we consider the vector space
$$
V_{\Lambda}=V_{\Lambda,\wg}=\{\v\in\wg^\bot\;\vert\; \forall j\not\in\Lambda,\; {\rm v}_j=0\}.
$$
Thus $V_\emptyset=\{0\}$ and $\dim V_\Lambda=\vert\Lambda\vert -1$ if $\Lambda\neq\emptyset$. We clearly have $V_{\Lambda_1}\cap V_{\Lambda_2}=V_{\Lambda_1\cap\Lambda_2}$.
\begin{remark}
\label{glue}
Let $\Lambda_1$, $\Lambda_2$ be non-empty subsets of $\{1,\ldots,r\}$. If $\Lambda_1\cap\Lambda_2\not=\emptyset$ 
we have
$$
V_{\Lambda_1}+V_{\Lambda_2}=V_{\Lambda_1\cup\Lambda_2}
$$
while $V_{\Lambda_1}+V_{\Lambda_2}$ is a subspace of $V_{\Lambda_1\cup\Lambda_2}$ of codimension $1$ if $\Lambda_1\cap\Lambda_2=\emptyset$.
\end{remark}
\DIM
Let us assume $\Lambda_1\cap\Lambda_2\not=\emptyset$. The displayed formula follows from the trivial inclusions $V_{\Lambda_i}\subseteq V_{\Lambda_1\cup\Lambda_2}$ and from the equality of dimensions:
\begin{align*}
\dim(V_{\Lambda_1}+V_{\Lambda_2})
&=\dim(V_{\Lambda_1})+\dim(V_{\Lambda_2})-\dim(V_{\Lambda_1}\cap V_{\Lambda_2})\\
&=\dim(V_{\Lambda_1})+\dim(V_{\Lambda_2})-\dim(V_{\Lambda_1\cap\Lambda_2})\\
&=(\vert\Lambda_1\vert-1)+(\vert\Lambda_2\vert-1)-(\vert\Lambda_1\cap\Lambda_2\vert-1)\\
&=\vert\Lambda_1\cup\Lambda_2\vert-1=\dim(V_{\Lambda_1\cup\Lambda_2}).
\end{align*}
The last assertion follows similarly.
\CVD

In order to generalise this simple remark, we introduce the following definition.
\begin{definition}
\label{connesso}
Let $\Gamma=\{\Lambda_1,\ldots,\Lambda_s\}$ be a collection of subsets of $\{1,\ldots,r\}$. We say that a subset $C$ of $\{1,\ldots,r\}$ is a connected component of $\Gamma$, if, after renumbering $\Lambda_1,\ldots,\Lambda_s$ if necessary, there exists an integer $k$ with $1\leq k\leq s$ such that
\begin{itemize}
\item[i)] $C=\Lambda_1\cup\cdots\cup\Lambda_k$;
\item[ii)] for $j=1,\ldots,k-1$ we have $\Lambda_{j+1}\cap(\Lambda_1\cup\cdots\cup\Lambda_j)\neq\emptyset$;
\item[iii)] for $j=k+1,\ldots,s$ we have $C\cap\Lambda_j=\emptyset$.
\end{itemize}
We say that $\Gamma$ is connected, if it has only one connected component.
\end{definition}

We may also rephrase this definition,  as follows: consider first the graph   on $\{1,\ldots ,s\}$ defined by joining $i, j$ if and only if $\Lambda_i\cap\Lambda_j\neq\emptyset$. Then a connected component in our sense is a union $C=\bigcup_{i\in U}\Lambda_i$, where $U$ is a connected component,  in the usual sense,  of the graph just defined.

\medskip

As an example, $\Gamma=\{\{1,2\},\{3,4,5\},\{2,5\},\{6,7,8\},\{7,8,9\}\}$ has two connected components, 
$\{1,2,3,4,5\}$ and $\{6,7,8,9\}$.\\

By Remark~\ref{glue} we easily see that:
\begin{remark}
\label{increase}
~
{\rm Let $\Gamma=\{\Lambda_1,\ldots,\Lambda_s\}$ be a collection of subsets of $\{1,\ldots,r\}$ and let $C_1,\ldots,C_p$ be the connected components of $\Gamma$. Then, 
\begin{itemize}
\item[i)] $V_{\Lambda_1}+\cdots+V_{\Lambda_s}=V_{C_1}+\cdots+V_{C_p}$.
\item[ii)] $\dim(V_{\Lambda_1}+\cdots+V_{\Lambda_s})=\vert\Lambda_1\cup\cdots\cdot\cup\Lambda_s\vert - p$.
\end{itemize}
}
\end{remark}

The following definition is crucial for our purposes.
\begin{definition}
Let $V$ be a $\Qb$-vector space, $v_1,\ldots,v_s$ be $s\geq 2$ vectors of $V$. Let $a_1v_1+\cdots+a_sv_s=0$ be a non-trivial linear relation. We say that this relation is {\it minimal} if there are no non-trivial relations $\sum b_i v_i=0$ over a proper non-empty subset of $\{1,\ldots,s\}$. 
\end{definition}
We remark that the relation $a_1v_1+\cdots+a_sv_s=0$ is minimal if and only if $a_1,\ldots,a_s\in\Qb^*$ and $\dim\langle v_1,\ldots,v_s\rangle=s-1$.

We also remark that, given $v_1,\ldots,v_s\in V$ linearly dependent  and not all zero, there exists a subset $\Lambda\subseteq\{1,\ldots,r\}$ such that $\{v_i\}_{i\in\Lambda}$ satisfy a minimal linear relation.\\

We now agree on some conventions which will be followed in the rest of this section and in the next section. 

We fix as above $r\geq 2$ rational functions $f_1,\ldots,f_r\in\F\backslash\{0\}$ and a non-constant $t\in\F$. We define $S_0$ as the finite set consisting of all zeros and poles of $f_1,\ldots,f_r,\d t$. We fix a point $Q\in\Cu(\Qb)\backslash S_0$ and we define $S=S_0\cup\{Q\}$. We choose once and for all  a positive integer $\Delta$, a real $c$ and $g,g_0,g_1,\ldots,g_{\Delta-1}\in L(\infty Q)$ depending only on $Q$ and $t$ and satisfying the statement of Lemma~\ref{good-basis}. The implicit constants in the big-$O$ below will depend only on these data.\\

The next lemma is the main tool in our construction. 

\begin{lemma}
\label{Wronskiano}
Let $\Lambda\subseteq\{1,\ldots,r\}$ be of cardinality $\geq2$ and $(M'_i)_{i\in\Lambda}$ be positive integers with maximum $M'$. Define
$$
d_\Lambda=-\deg\div(f_i)_{i\in\Lambda},\qquad \Theta=\max\Big(1,\sum_{i\in\Lambda} M'_i-(M'+nd_\Lambda)\Big).
$$
Let also $\Pt\in\Cu(\Qb)\backslash S$, $n\in\N$ and put 
$$
\wg=(f_1(\Pt)^n,\ldots,f_r(\Pt)^n).
$$ 
Let finally $\{A_i\}_{i\in\Lambda}\subset\F$ not all zero, such that $A_i\in L(M'_iQ)$. Let us assume that $\{A_if_i^n\}_{i\in\Lambda}$ satisfy a minimal linear relation. Then there exists a basis of algebraic vectors $\v_1,\ldots,\v_{\vert\Lambda\vert-1}$ of $V_{\Lambda,\wg}$ which satisfies
$$
h(\v_j)\leq M' h(\Pt)+O\left(\Theta (h(\Pt)+1)+(n+M')(1+h(\Pt)^{1/2})+\max h(A_i)\right).
$$
\end{lemma}

\DIM In the proof we use a Wronskian argument. Let us first recall some basic facts on it. The derivative $\frac{d}{d t}$ on $\Qb(t)$ can be uniquely extended to $\F$. For $F\in\F$ and $R\in\Cu(\Qb)$ we have 
\begin{equation}
\label{local-der}
\begin{cases}
\ord_R(\d F/\d t)=\ord_R(F)-1-\ord_R(\d t),& \hbox { if }\ord_R(F)\neq0 ;\\
\ord_R(\d F/\d t)\geq -\ord_R(\d t),& \hbox { if }\ord_R(F)=0 .
\end{cases}
\end{equation}
The (normalized) {\it Wronskian} of $F_1,\ldots,F_k\in\F$ with respect to $t$ is the determinant 
$$
W(F_1,\ldots,F_k)=\det\left(\frac{1}{j!}\frac{\d^jF_i}{\d t^j}\right)_{\newatop{i=1,\ldots,k}{j=0,\ldots,k-1}}.
$$
It is well known that $W = 0$ if and only if the $F_i$'s are linearly dependent over $\Qb$. 

\medskip

Let us now go on with the proof of our lemma. We may assume $\Lambda=\{1,\ldots,s\}$. For short we put $F_i=A_if_i^n$ and 
$W_i=W(F_1,\ldots,F_{i-1},F_{i+1},\ldots,F_s)$ for $i=1,\ldots,s$. 

For later reference, we remark that $F_1,\ldots,F_s$ are $S$-integers (as elements of the function field $\F$): indeed the $\div(f_i)$ are supported in $S$ and the $A_i$ are also $S$-integers since $A_i\in L(M'Q)$.  Moreover the zeros and the poles of $dt$ are in $S$ as well. Thus $\d^l F_i/\d t^l$ are $S$-integers (\cf \eqref{local-der}), and so also~$W_s$.

\smallskip

By assumption, we have a minimal linear relation $a_1F_1+\cdots+a_sF_s=0$. Thus $F_1,\ldots,F_{s-1}$ are linearly independent and  $a_1,\ldots,a_s\neq0$. This proves that $W_s\neq0$. 

Let $1\leq i\leq s-1$. Since $a_i\neq0$ we can replace $F_i$ by 
$$
-\frac{a_1}{a_i}F_1-\cdots-\frac{a_{i-1}}{a_i}F_{i-1}-\frac{a_{i+1}}{a_i}F_{i+1}-\cdots-\frac{a_s}{a_i}F_s
$$
in $W_s$. This shows that $W_s=\pm (a_s/a_i)W_i$. 

\smallskip

We want to obtain a suitable upper bound for $m_0:=\ord_{\Pt}(W_s)$. For this, we shall  use the fact that $W_s$ has already a big multiplicity at the zeros of $F_i$, since these functions are essentially $n$-th powers.

For $i=1,\ldots,s$ and for $l=0,\ldots,s-1$ we have (\cf \eqref{local-der}) 
$$
\ord_R\left(\frac{\d^l F_i}{\d t^l}\right)\geq \ord_R(F_i)-l(1+\ord_R(\d t))=\ord_R(A_i)+n\, \ord_R(f_i)+O(1)
$$
for any $R\in\Cu(\Qb)$. Thus, for $i=1,\ldots,s$, 
\begin{align*}
\ord_R(W_s)=\ord_R(W_i)
&\geq \sum_{j\neq i}\Big(\ord_R(A_j)+n\, \ord_R(f_j)+O(1)\Big)\\
&=\sum_{j=1}^sn\, \ord_R(f_j)-\lambda_{R,i}+O(1).
\end{align*}
where $\lambda_{R,i}:=n\, \ord_R(f_i)-\sum_{j\neq i}\ord_R(A_j)$. We deduce:
\begin{align*}
\ord_R(W_s)
&\geq \max_{i=1,\ldots,s}\left(\sum_{j=1}^sn\, \ord_R(f_j)-\lambda_{R,i}+O(1)\right)\\
&=\sum_{j=1}^sn\, \ord_R(f_j)-\lambda_R+O(1)
\end{align*}
where we have defined, for $R\in S$,
$$
\lambda_R:=\min_{i=1,\ldots,s}\lambda_{R,i}=\min_{i=1,\ldots,s}\Big\{n\, \ord_R(f_i)-\sum_{j\neq i}\ord_R(A_j)\Big\}.
$$
We shall use this inequality,  for $R\in S$,  in the functional  (\ie in $\F$) product formula, applied to $W_s$, namely the formula $\sum_R\ord_R(W_s)=0$.

In this formula, for $R=\Pt\not\in S$ we find the quantity  $\ord_{\Pt}(W_s)$ that we have to estimate, whereas  for $R$ outside $S\cup\{\Pt\}$ we use the trivial bound $\ord_R(W_s)\geq0$. Also, since the $\div(f_i)$ are supported in $S_0\subseteq S$, we have\footnote{We use the fact that $S$ contains not only the poles but also the zeros of the $f_i$'s.} $\sum_{R\in S}\ord_R(f_i)=\deg(\div(f_i))=0$. Moreover $\Pt\not\in S$. Thus 
$$
0=\sum_{all\ R}\ord_{R}(W_s)
\geq\ord_{\Pt}(W_s)+\sum_{R\in S}\ord_R(W_s)
\geq\ord_{\Pt}(W_s)-\sum_{R\in S}\lambda_R+O(1).
$$
We now recall that $A_i\in L(M'_iQ)$ and that the $f_i$ are supported in $S_0$. Thus 
$$
\lambda_Q=\min_j\Big\{-\sum_{i\neq j}\ord_Q(A_i)\Big\}\leq \min_j\sum_{i\neq j} M'_i=  \sum_{i=1}^sM'_i-M'
$$
and
$$
\sum_{R\in S_0}\lambda_R
\leq \sum_{R\in S_0}\min_i\{n\, \ord_R(f_i)\}
=n\sum_{R}\min_i\{\ord_R(f_i)\}=
-nd_\Lambda.
$$
Collecting together these last three inequalities, we get 
the following sought upper bound for $m_0=\ord_{\Pt}(W_s)$:
\begin{equation}
\label{m0}
m_0\leq\sum_{R\in S}\lambda_R+O(1)
\leq\sum_{i\in\Lambda} M'_i-(M'+nd_\Lambda)+O(1)=\Theta+O(1).
\end{equation}
\medskip

We can now construct the desired basis $\v_1,\ldots,\v_{\vert\Lambda\vert-1}$ of $V_{\Lambda}=V_{\Lambda,\wg}$. For a non-negative integer $\rho$ we put as before
$$
\delta_\rho=\frac{1}{\rho!}\frac{\d^{\rho}}{\d t^\rho}
$$
Given a vector of non-negative integers $\rhog=(\rho_1,\ldots,\rho_{s-1})$ we let
$$
W_{\rhog}=\det(\delta_{\rho_j}F_i)_{i,j=1,\ldots,s-1}.
$$
Thus $W_s=W_{(0,1,\ldots,s-1)}$. It is also easily seen that
\begin{equation}
\label{deco}
\delta_{m_0}W_s\in \sum_{\vert\rhog\vert=e} \Z W_{\rhog}
\end{equation}
where $\vert\rhog\vert=\rho_1+\cdots+\rho_{s-1}$ and $e=m_0+1+\cdots+(s-2)$. By~(\ref{m0}) we have 
\begin{equation}
\label{multiplicity}
e\leq\Theta+O(1).
\end{equation}

Since $\Pt\not\in S$ and the zeros of $dt$ are in $S$, we have $\ord_{\Pt}(dt)=0$. By definition of $m_0$ and by~\eqref{local-der} we have $\ord_{\Pt}(\delta_{m_0}W_s)=0$. Again by~(\ref{local-der}) and since $F_i$ are $S$-integers, for all $\rhog$ we have $\ord_{\Pt}(\delta_{m_0}W_{\rhog})\geq 0$. By~(\ref{deco}), this implies that there exists $\rhog'$ with $\vert\rhog'\vert=e$ such that $\ord_{\Pt}(W_{\rhog'})=0$. For $i=1,\ldots,s$ and $j=1,\ldots,s-1$, let 
\begin{equation}
\label{B}
B_{ij}=a_i f_i^{-n}\delta_{\rho'_j}F_i.
\end{equation}

Recall that: $\Pt\not\in S$, the $f_i$'s are supported in $S$, the $F_i$ have all their poles in $S$, the zeros of $dt$ are in $S$. By~\eqref{local-der} we see that $\ord_{\Pt}(B_{ij})\geq0$. Thus $B_{ij}(\Pt)\in\Qb$. 

Since $a_1F_1+\cdots+a_sF_s=0$, we have $B_{1j} f_1^n+\cdots+B_{sj} f_s^n=0$ for $j=1,\ldots,s-1$. Thus, for $i=1,\ldots,s$, 
$$
\v_j=(B_{1j}(\Pt),\ldots,B_{sj}(\Pt),0,\ldots,0)\in V_{\{1,\ldots,s\}}.
$$

\smallskip

The important fact that we have achieved so far is that {\sl since $\Pt$ is not a zero of $W_{\rhog'}$ the vectors $\v_1,\ldots,\v_{s-1}$ are linearly independent and form a basis of $V_{\{1,\ldots,s\}}$.}

\medskip

By Lemma~\ref{height-stuff}, $h(F_i)=O(n+M_i'+h(A_i))$. In order to deduce an upper bound for the height of $B_{ij}$ we still need a bound for the height of the coefficients $a_1,\ldots,a_s$ of the minimal linear relation $a_1F_1+\cdots+a_sF_s=0$. Obviously, we may assume $a_s=-1$. We differentiate the relation up to order $s-2$. Since the rational functions $F_1,\ldots,F_{s-1}$ are linearly independent over $\Qb$, their Wronskian is not zero and so we can solve the resulting system for $a_1,\ldots,a_{s-1}$. Using Cram\'er's Rule and $h(F_i)=O(n+M_i'+h(A_i))$  together with Lemma \ref{height-stuff} especially (i) and (ii) we find without difficulty that $h(a_i)=O(n+M'+\max h(A_i))$.

Now to simplify the notation we write $\v_j$ as $(B_{1}(\Pt),\ldots,B_{s}(\Pt),0,\ldots,0)$ for $B_i=a_if_i^{-n}\delta_\rho(A_if_i^n)$ as in (\ref{B}); here $\rho \leq e \leq \Theta+O(1)$ by (\ref{multiplicity}). As before
$$f_i^{-n}\delta_\rho(A_if_i^n)=f_i^{-n}\sum_{l+m=\rho}\delta_l(A_i)\delta_m(f_i^n)$$
which by~\eqref{francesco} is
$$\sum_{l+m=\rho}\delta_l(A_i)f_i^{-m}\sum_{\bf a}C({\bf a})\delta_0(f_i)^{a_0}\delta_1(f_i)^{a_1}\cdots\delta_m(f_i)^{a_m}$$
(note the changed power of $f_i$). As in Lemma \ref{good-basis}(ii) we write
$$A_i=\sum_{j,k}\alpha_{ijk}g_jg^k~~~~(i=1,\ldots,s)$$
and it suffices here to take $k \leq M'_i/\Delta\leq M'/\Delta$. Again~\eqref{francesco} for $\delta_q(g^k)$ gives 
$$f_i^{-n}\delta_\rho(A_if_i^n)=\sum_{l+m=\rho}f_i^{-m}\sum_{p+q=l}\sum_{j,k}\alpha_{ijk}\delta_p(g_j)g^{k-q}E_{qm},$$
where
\begin{multline}
\label{star}
E_{qm}=\\ \sum_{\bf b}C({\bf b})\delta_0(g)^{b_0}\delta_1(g)^{b_1}\cdots\delta_q(g)^{b_q}\sum_{\bf a}C({\bf a})\delta_0(f_i)^{a_0}\delta_1(f_i)^{a_1}\cdots\delta_m(f_i)^{a_m}.
\end{multline}

We evaluate all this at $\Pt$ and we want a main term $M'h(\Pt)$ ``uniformly in $i$", so that when we bundle the coordinates into the vector we don't get $sM'h(\Pt)$.

In fact the terms $g^{k-q}$ with $k\geq q$ already give
$$
h\Big(g^{k-q}(P)\Big)=(k-q)h(g(\Pt))\leq {M' \over \Delta}h(g(\Pt))\leq {M' \over \Delta}\Big(\Delta h(\Pt)+O(1+\sqrt{h(\Pt)})\Big)
$$
by the upper bound in Lemma \ref{HeightMachine} with just two functions. The right-hand side is
$$
M'h(\Pt)+O\big(M'(1+\sqrt{h(\Pt}))\big)
$$
so we already have the main term, clearly uniformly. Thus the rest had better be small. The point here is 
\begin{align*}
b_0+b_1+\cdots+b_q&= b_1+\cdots+qb_q= q \leq l \leq \rho \leq \Theta+O(1),\\
a_0+a_1+\cdots+a_m&= a_1+\cdots+ma_m= m \leq\rho\leq \Theta+O(1).
\end{align*}
So if $k<q$ in $g^{k-q}$ then $|k-q|\leq q$ so we get
$$
h\Big(g^{k-q}(P)\Big)=|k-q|h(g(\Pt))\leq qh(g(\Pt))=O\left(\Theta(h(\Pt)+1)\right).
$$
And by Lemma \ref{heisenstein} we get for the $\delta$-terms in (\ref{star}), as well as the $f_i$-term, a height of order at most
$$
(q+m)(h(\Pt)+1)\leq \rho(h(\Pt)+1)\leq (\Theta+O(1))(h(\Pt)+1),
$$
also uniformly. The $C$-terms contribute logarithmically to order at most
$$
k+q+n+m\leq k+n+\rho\leq M'+n+\Theta+O(1).
$$
So this deals uniformly with the $E_{qm}$. As $p\leq\rho$ the $\delta_p(g_j) ~(j=0,1,\ldots,\Delta)$ give nothing new, and by Lemma \ref{good-basis}(ii) the $\alpha$-terms contribute $O(M'+\max h(A_i))$. 

Collecting everything up, we get
$$
h(\v_j)\leq M'h(\Pt)+O\left(\Theta (h(\Pt)+1)+n+\max h(A_i)+M'(1+h(\Pt)^{1/2})\right)
$$
for $j=1,\ldots,s-1$, slightly better than required. 
\CVD

\section{Proof of Theorem~\ref{main}.}
\label{strategy}

The next theorem is a refined version of Theorem~\ref{main}.

\begin{theorem}
\label{main-explicit}
Let $r\geq 2$ and $f_1,\ldots,f_r\in\F$ be non-zero rational functions such that $f_i/f_j$ is non-constant for some $i\neq j$. Let $d=-\deg\div(f_1,\ldots,f_r)$, $K>0$ be sufficiently large with respect to $f_1,\ldots,f_r$ and $\alphag=(\alpha_1:\cdots:\alpha_r)\in\P^{r-1}(\Qb)$. Consider, for a natural number $n$, a solution $P\in \Cu(\Qb)$ of the equation
$$
\alpha_1 f_1(\Pt)^n+\cdots + \alpha_r f_r(\Pt)^n =0.
$$
Then, if $n\geq K$ and if there are no proper vanishing subsums, we have 
\begin{equation}
\label{main-inequa}
h(\Pt)\leq \left(\frac{r-1}{d}+O(1/K)\right)\frac{h(\alphag)}{n}+O(K^2)
\end{equation}
 where the implicit constant in the big-$O$ depends only on $f_1,\ldots,f_r$.
\end{theorem}

\medskip

Theorem~\ref{main} easily follows from Theorem~\ref{main-explicit}. Indeed, choosing $K$ sufficiently large, we have $(r-1)/d+O(1/K)\leq (r-1)/d +1\leq r$.\\

\begin{remark}
\label{comments}~\\
{\rm 
i) In the proof we shall show:
\begin{equation}
\label{main-inequa2}
\frac{h(\alphag)}{n}
\geq\frac{d h(\Pt)}{r-1}+O\left(\frac{1}{K}h(\Pt)+h(\Pt)^{1/2}+K\right).
\end{equation}
which immediately implies~\eqref{main-inequa}, since either $h(\Pt)\leq K^2$ or $h(\Pt)^{1/2}+K\leq 2h(\Pt)/K$. Note also that we will not use the assumption $f_i/f_j\neq {\rm constant}$ in the proof of~\eqref{main-inequa2}. This assumption is equivalent to $d\neq0$ and~\eqref{main-inequa2} is trivially satisfied if $d=0$.\\

\noindent ii) Let $\betag=(f_1(\Pt):\cdots:f_r(\Pt))\in\P^{r-1}(\Qb)$. By Weil's Height Machine Lemma~\ref{HeightMachine} (and since $d>0$), inequality~\eqref{main-inequa2} is equivalent to  
\begin{equation}
\label{main-inequa3}
\frac{h(\alphag)}{n}
\geq\frac{h(\betag)}{r-1}+O\left(\frac{1}{K}h(\betag)+h(\betag)^{1/2}+K\right)
\end{equation}
(which in turn implies $h(\betag)\leq (r-1+O(1/K))\frac{h(\alphag)}{n}+O(K^2)$).\\

\noindent iii) A standard application of Siegel's lemma to the linear equation 
\begin{equation}
\label{linear}
\alpha_1\beta_1^n+\cdots+\alpha_r\beta_r^n=0
\end{equation}
in the unknowns $\alpha_1,\ldots,\alpha_r$, shows that there exists a solution with
$$
\frac{h(\alphag)}{n}\leq \frac{h(\betag)}{r-1}+O(1).
$$
Thus inequality~\eqref{main-inequa3} is sharp. More precisely, Theorem~\ref{main-explicit} gives a lower bound for the first minimum  (with respect to the height) of the linear equation~\eqref{linear}, and hence shows that the successive minima are close to each other for large $n$.\\

\noindent iv) 
 If we ask that $f_i/f_j$ is non-constant for {\sl all} $i\neq j$, the assumption on vanishing subsums can be easily removed (by induction on $r$).\\

\noindent v) We finally remark that  we have a result  even if $d=0$ (\ie $f_i/f_j$ constant for all $i$, $j$), but now the lower bound $n\geq K$ must depend also on $\alphag$. Indeed, if $f_i=c_i f_1$ with $c_i$ constants, our equation becomes $(\alpha_1 c_1^n+\cdots + \alpha_r c_r^n)f_1(P)^n=0$. By the Skolem-Mahler-Lech Theorem, if $f_1(P)\neq0$ then $n$ is bounded by a constant depending on $\alpha_1,\ldots,\alpha_r$ and on $c_1,\ldots,c_r$.}
\end{remark}

{\smallskip\noindent{\bf Strategy of the proof of Theorem~\ref{main-explicit}.}\quad}
Let $f_1,\ldots,f_r\in\F$ be as in the statement of the theorem. We recall that we have chosen a non-constant rational function $t\in\F$ and that $S$ is the finite set consisting of all zeros and poles of $f_1,\ldots,f_r,dt$ and of an extra point $Q$ (which is neither a zero nor a pole of $f_1,\ldots,f_r,dt$).

We fix algebraic numbers $\alpha_1,\ldots,\alpha_r$, not all zero. In order to prove~\eqref{main-inequa} we may suppose that $\Pt$ does not lie in any prescribed finite set of points. We thus choose $\Pt\in\Cu(\Qb)\backslash S$ satisfying our equation 
$$
\alpha_1 f_1(\Pt)^n+\cdots + \alpha_r f_r(\Pt)^n =0
$$ 
for some $n\geq K$. We shall also assume $h(\Pt)\geq1$.

\medskip

Put now $\wg=(f_1^n(\Pt),\ldots,f_r^n(\Pt))$. Thus $\alphag\in\wg^\bot$, the orthogonal space of $\wg$. 
Our strategy is the following:  

We shall first construct a basis of function-vectors (with controlled heights) for the orthogonal of the vector $(f_1^n,\ldots ,f_r^n)\in\F^r$. Then we shall specialise at $P$, in order to obtain a basis of $\wg^\bot$, again with controlled heights. All of this shall involve an induction, necessary to take into account certain unexpected  linear relations, \ie relations with certain special properties in addition to those imposed by the construction. 

At this stage we shall get a new basis, on replacing   one of the vectors of the previous  basis  with $\alphag$. By well-known facts,  $\wg$ and $\wg^\bot$ have the same height. Lemma~\ref{HeightMachine} gives a lower bound for the height of $\wg$. The height of $\wg^\bot$ is bounded from above by the sum of the heights of the vectors of our new basis. Comparing these bounds, we shall get the desired conclusion.\\

{\smallskip\noindent{\bf First step of the inductive proof.}\quad}

Let $N_1$ be the minimum of the set of integers $m\geq0$ such that there exist a non-empty $\Lambda\subseteq\{1,\ldots,r\}$ and 
rational functions $A_i\in\F$ ($i\in\Lambda$) not all zero, satisfying
\begin{equation}
\label{N1}
\begin{cases}
A_i\in L(mQ),& \hbox{ for }Êi\in\Lambda;\\ 
h(A_i)\leq nK,& \hbox{ for } i\in\Lambda;\\
(A_if_i^n)_{i\in\Lambda} \hbox{ are linearly dependent over $\Qb$.}
\end{cases}
\end{equation}\\

\begin{fact}
\begin{equation}
\label{N1Siegel}
(r-1)N_1\leq nd + O(n/K).
\end{equation}
\end{fact}
\DIM
We provide an upper bound for $N_1$ using Lemma~\ref{siegel-substitute},  as we are going to illustrate. Let $0<\e<1/2$ and for this argument define $N$ as the smallest integer such that 
\begin{equation}
\label{Ndef}
(r-1-\e)N\geq (1+\e)nd + c_0.
\end{equation}
Then $N=O(n)$ and, with $M_1=\cdots=M_r=N$,  the Dirichlet exponent $\varrho$ of Lemma~\ref{siegel-substitute} satisfies
$$
\varrho=\frac{N+dn}{(r-1)N-dn-c_0}\leq\frac{1}{\e}.
$$
Thus, by that lemma, there exist $A_1,\ldots, A_r\in L(N Q)$ not all zero such that 
\begin{equation}
\label{*2}
A_1f_1^n+\ldots+A_rf_r^n=0
\end{equation}
and 
$$
h(A_i)= O(n/\e+N/\e)=O(n/\e).
$$
Choosing $\e=c/K$, where $c$ is a sufficiently large constant to kill the implicit constant in the last $O()$, we see that there exists a non-trivial solution of~(\ref{*2}) with $A_i\in L(NQ)$ and $h(A_i)\leq nK$. Not all $A_1,\ldots,A_r$ are zero,  and we see that the non-zero ones among $A_1f_1^n,\ldots,A_rf_r^n$ sum up to zero, and so they are linearly dependent, as required by~(\ref{N1}) (on choosing $\Lambda$ simply as the set of $i$ such that $A_i\neq0$). 

This shows that $N_1\leq N$. Since $N$ is the smallest integer satisfying~(\ref{Ndef}) and since $\e=c/K$, we have $(r-1)N= nd +O(n/K)$. Thus~\eqref{N1Siegel} holds.
\CVD

Among all subsets $\Lambda$ which realize the minimum defining $N_1$ in~(\ref{N1}), we choose a subset which is minimal. We denote by $\Lambda_1$ such a set, by $l_1=\vert\Lambda_1\vert$ its cardinality (necessarily $l_1\ge 2$) and  by $\{A^{(1)}_i\}_{i\in\Lambda_1}$ the corresponding rational functions. This implies in particular that $\{A^{(1)}_i f_i^n\}_{i\in\Lambda_1}$ satisfy a minimal linear relation.\\

To go ahead we want to apply Lemma~\ref{Wronskiano} to find a suitable basis of $V_{\Lambda_1}$. In that lemma, let us put  $M'_i=N_1$ for $i\in\Lambda_1$, so $M'=\max M'_i=N_1=O(n)$ and 
$$
\Theta=\max\Big(1,\sum_{i\in\Lambda_1}M'_i-(M'+nd_1)\Big)=\max\Big(1,(l_1-1)N_1 -nd_1\Big)
$$ 
with $d_1=d_{\Lambda_1}=-\deg\div(f_i)_{i\in\Lambda_1}$. 

\begin{fact}
$\Theta=O(n/K)$. 
\end{fact}
\DIM To prove this, we use again Lemma~\ref{siegel-substitute}, this time on the $f_i~(i\in\Lambda_1)$, with $M_i=N_1-1$ for $i\in\Lambda_1$. The Dirichlet exponent $\varrho$ is then
$$
\varrho=\frac{N_1-1+nd_1}{(l_1-1)(N_1-1)-nd_1-c_0}.
$$
By Lemma~\ref{siegel-substitute}, there exist rational functions $B_i$ not all zero such that $B_i\in L((N_1-1) Q)$ for $i\in\Lambda_1$,  $\sum_{i\in\Lambda_1} B_if_i^n=0$ and 
$$
h(B_i)= O((\varrho+1) n).
$$

By the minimality of $N_1$, we cannot have $\max_i h(B_i)\leq nK$. Thus $\varrho\geq K/c$, where $c$ is a sufficiently large constant to kill the implicit constant in the last $O()$. This implies
$$
(l_1-1)N_1-nd_1\leq \frac{c}{K}(N_1+nd_1)+O(1)=O(n/K)
$$
as required. This concludes the proof of this fact.
\CVD

Using Lemma~\ref{Wronskiano} and the inequalities
$$
h(\Pt) \geq 1,\quad h(A_i)\leq nK,\quad M'_i=N_1=O(n),\quad \Theta=O(n/K),
$$
we find a basis $\v^{(1)}_1,\ldots,\v^{(1)}_{l_1-1}$ of $V_{\Lambda_1}$ satisfying
\begin{equation}
\label{basisV1}
\begin{aligned}
h(\v^{(1)}_i)
&\leq M' h(\Pt)\\
&\hskip 0.25cm+O\left(\Theta (h(\Pt)+1)+(n+M')(1+h(\Pt)^{1/2})+\max h(A_i)\right)\\
&= N_1h(\Pt)+O\left(\frac{n}{K}h(\Pt)+nh(\Pt)^{1/2}+nK\right).
\end{aligned}
\end{equation}\\

{\smallskip\noindent{\bf Proof of Theorem~\ref{main-explicit} in a test-run case.}\quad
Let as now assume $\Lambda_1=\{1,\ldots,r\}$, which is in essence the generic case. In this case we shall obtain directly the desired conclusion of Theorem~\ref{main-explicit}, as we now show.}\\

Recall that $\wg=(f_1^n(\Pt),\ldots,f_r^n(\Pt))$. Since $\alphag$ is a non-zero vector in $\wg^\bot$ we may assume (reordering $\v^{(1)}_1,\ldots,\v^{(1)}_{r-1}$ if necessary) that 
$$
\alphag,\v^{(1)}_1,\ldots,\v^{(1)}_{r-2}
$$
is a basis of $\wg^\bot$. Let us denote by $h_2$ the logarithmic euclidean height (defined  on choosing the $L_2$-norm at the infinite places). By well-known facts on the height of subspaces (see~\cite{Bo-Va} and~\cite{Sc2}) and by the previous upper bounds for the height of these vectors. 
\begin{align*}
h(\wg)
&\leq  h_2(\wg)=h_2(\wg^\bot)\leq h_2(\alphag)+\sum_{i=1}^{r-2}h_2(\v^{(1)}_i)+\log(r-1)\\
&\leq (r-2)N_1h(\Pt)+h(\alphag)+O\left(\frac{n}{K}h(\Pt)+nh(\Pt)^{1/2}+nK\right).
\end{align*}
Moreover, by the functorial lower bound for the height Lemma~\ref{HeightMachine}, we have:
$$
h(\wg)\geq \left(h(\Pt)+O(h(\Pt)^{1/2})\right)nd.
$$
Thus
$$
0\leq\lambda h(\Pt)+ \frac{h(\alphag)}{n}+O\left(\frac{1}{K}h(\Pt)+h(\Pt)^{1/2}+K\right)
$$
with
$$
\lambda=(r-2)N_1/n-d.
$$
By~(\ref{N1Siegel}) we have 
$$
\lambda\leq (r-2)\frac{d}{r-1}-d+O(1/K)=-\frac{d}{r-1}+O(1/K).
$$
Inequality~\eqref{main-inequa2} follows.\\

{\smallskip\noindent{\bf Inductive construction.}\quad
The obstacle in the approach of the test-run case is that $\Lambda_1$ may be smaller than $\{1,\ldots ,r\}$. If this happens, we can somewhat take advantage of the fact that we have an `unexpected' dependence relation. To exploit this, let us sketch how we intend to argue by an induction procedure. The following claim resumes the inductive construction we shall do in the next section.}

\begin{claim}
\label{claim}
There exists an integer $s$ with $1\leq s\leq r$, positive integers $N_1,\ldots\, N_s$ and non-empty subsets $\Lambda_1,\Lambda_2,\ldots,\Lambda_s$ of $\{1,\ldots,r\}$ of cardinalities $l_1,l_2\ldots,l_s$ satisfying: 
\begin{itemize}
\item[i)] $N_1\leq N_2\leq\cdots\leq N_s$.
\item[ii)] For $j=2,\ldots,s$, the set $\Lambda_j$ is contained in no connected component of $\{\Lambda_1,\ldots,\Lambda_{j-1}\}$.
\item[iii)] The collection $\{\Lambda_1,\ldots,\Lambda_s\}$ is connected and its union is the full set $\{1,\ldots,r\}$.
\item[iv)] Let $t_1=\dim(V_{\Lambda_1})$ and, for $j=2,\ldots,s$ let
$$
t_j=\dim(V_{\Lambda_1}+\cdots+V_{\Lambda_j})-\dim(V_{\Lambda_1}+\cdots+V_{\Lambda_{j-1}}).
$$
Then for $j=1,\ldots,s$ we have
$$
t_1N_1+\cdots+t_{j-1}N_{j-1}+\Big(r-1-\sum_{i=1}^{j-1}t_i\Big)N_j\leq nd+ O(n/K).
$$
\item[v)] For $j=1,\ldots,s$, there exists a basis $\v^{(j)}_1,\ldots,\v^{(j)}_{l_j-1}$ of $V_{\Lambda_j}$ satisfying
$$
\max_ih(\v^{(j)}_i)\leq N_jh(\Pt)+O\left(\frac{n}{K}h(\Pt)+nh(\Pt)^{1/2}+nK\right).
$$
\end{itemize}
\end{claim}~\\

We shall explain how to perform this construction in the next section, proving the claim. For the moment, we pause to show how this claim allows us to conclude the proof of Theorem~\ref{main-explicit}.\\

We state at once a general elementary lemma.

\begin{lemma}
\label{convex}
Let $\tau$, $\rho$, $a_1,\ldots,a_s$ be positive real numbers such that $a_s\geq 1$ and $a_1+\ldots+a_s=\rho$. Let also $x_1\leq\ldots\leq x_s$ be positive real numbers such that 
$$
\sum_{i=1}^s a_i x_i \leq \tau.
$$
Then 
$$
\sum_{j=1}^{s-1} a_j x_j+(a_s-1)x_s -\tau\leq -\tau/\rho.
$$
\end{lemma}
\DIM  Set $\sigma:=\sum_{i=1}^sa_ix_i$, so $\sigma\le \tau$ and also $\sigma\le \rho x_s$.   
Hence $\sigma-x_s\le \sigma(1-{1\over \rho})\le \tau(1-{1\over \rho})$, since  $\rho\ge a_s\ge 1$. 

Now, on subtracting $\tau$ from both sides we obtain $\sigma-x_s-\tau\le -{\tau\over \rho}$, as required.  
 \CVD
 
{\smallskip\noindent{\bf Deduction of Theorem~\ref{main-explicit} from Claim~\ref{claim}.}\quad}
We remark that
\begin{equation}
\label{note}
V_{\Lambda_1}+\cdots+V_{\Lambda_s}=\wg^\bot. \quad\hbox{In particular, } t_1+\cdots+t_s=r-1, 
\end{equation}
by Claim~\ref{claim} iii) and by Remark~\ref{increase} i). Thus, taking into account i) of Claim \ref{claim},  we may assume, after reordering for each $j$ the vectors $\v^{(j)}_1,\ldots,\v^{(j)}_{l_j-1}$ and possibly omitting some of them, that 
$$
\v^{(1)}_1,\ldots,\v^{(1)}_{t_1},\ldots,\v^{(s-1)}_1,\ldots,\v^{(s-1)}_{t_{s-1}},\v^{(s)}_1,\ldots,\v^{(s)}_{t_s}
$$
is a basis of $\wg^\bot$. 

We also remark that $t_s\geq 1$. Otherwise $\dim(V_{\Lambda_1}+\cdots+V_{\Lambda_{s-1}})=t_1+\cdots+t_{s-1}=r-1$ by~\eqref{note}. By Remark~\ref{increase} ii), this implies that $\{\Lambda_1,\ldots,\Lambda_{s-1}\}$ is connected and its union is the full set $\{1,\ldots,r\}$, which contradicts Claim~\ref{claim} ii) with $j=s$.

By assumption there are no proper  vanishing subsums in 
$\alpha_1 f_1(\Pt)^n+\cdots + \alpha_r f_r(\Pt)^n =0$. This implies that 
$\alphag\not\in V_{\Lambda_1}+\cdots+V_{\Lambda_{s-1}}$.
Thus we may assume that 
$$
\v^{(1)}_1,\ldots,\v^{(1)}_{t_1},\ldots,\v^{(s-1)}_1,\ldots,\v^{(s-1)}_{t_{s-1}}, \alphag,
\v^{(s)}_1,\ldots,\v^{(s)}_{t_s-1}
$$
is a basis of $\wg^\bot$. Arguing as we did before, we deduce that
\begin{equation}
\label{lambda}
0\leq\lambda h(\Pt)+\frac{h(\alphag)}{n}+O\left(\frac{1}{K}h(\Pt)+h(\Pt)^{1/2}+K\right)
\end{equation}
with
$$
\lambda=t_1\frac{N_1}{n}+\cdots+t_{s-1}\frac{N_{s-1}}{n}+(t_s-1)\frac{N_s}{n}-d.
$$

To go ahead, we recall that $t_1+\cdots+t_s=r-1$ (see~\eqref{note}). Thus, Claim~\ref{claim} iv) for $j=s$ reads
$$
t_1\frac{N_1}{n}+\cdots+t_{s-1}\frac{N_{s-1}}{n}+t_s\frac{N_s}{n}\leq d+O(1/K).
$$
We  apply the lemma~\ref{convex} with 
$$
a_j=t_j,\quad\rho=r-1,\quad x_j=N_j/n\quad{\rm and}\quad\tau=d+O(1/K).
$$
We find 
$$
\lambda=t_1\frac{N_1}{n}+\cdots+t_{s-1}\frac{N_{s-1}}{n}+(t_s-1)\frac{N_s}{n}-d\leq-d/(r-1)+O(1/K).
$$ 
Thus, by~\eqref{lambda},
$$
0\leq-\frac{d}{r-1}h(\Pt)+\frac{h(\alphag)}{n}+O\left(\frac{1}{K}h(\Pt)+h(\Pt)^{1/2}+K\right).
$$
Inequality~\eqref{main-inequa2} follows. This concludes the proof of Theorem~\ref{main-explicit}, assuming the truth of Claim~\ref{claim}.\\

\section{Proof of Claim~\ref{claim}}
\label{proofmain}
In this section, as promised, we detail our inductive process, verifying all the assertions of Claim~\ref{claim}.\\
We construct by induction an integer $s$ with $1\leq s\leq r$ and, for each $j=1,\ldots,s$, 
\begin{itemize}
\item[-] a positive integer $N_j$; 
\item[-] a subset $\Lambda_j$ of $\{1,\ldots,r\}$ of cardinality denoted $l_j:=|\Lambda_j|$; 
\item[-] a subset $J_j$ of $\Lambda_1\cup\cdots\cup\Lambda_{j-1}$  such that:
\begin{itemize}
\item[--] $\vert J_j\cap C\vert=1$ for each connected component of $\{\Lambda_1,\ldots,\Lambda_{j-1}\}$
\item[--] $J_j$ is disjoint from $\Lambda_1\cup\cdots \cup\Lambda_{j-2}\setminus J_{j-1}$ if $j\geq 2$.
\end{itemize}
\item[-] rational functions $\{A^{(j)}_i\}_{i\in\Lambda_j}$.
\item[-] a function $\varphi_j\colon\{1,\ldots,r\}\rightarrow\{1,\ldots,j\}$;
\end{itemize}

The role of the functions $\varphi_j$ shall appear along the discussion.\\

For $j=1$, we let $N_1$ be the minimum of the set of integers $N\geq0$ such that there exist $\Lambda\subseteq\{1,\ldots,r\}$ and 
rational functions $A_i\in\F$ ($i\in\Lambda$) not all zero, of height $\leq nK$, with $(A_if_i^n)_{i\in\Lambda}$ linearly dependent and such that
$$
A_i\in L(NQ),\hbox{ for }Êi\in\Lambda.
$$
Among all subsets $\Lambda$ which realize the minimum defining $N_1$, we choose a subset which is minimal. We denote $\Lambda_1$ such set, $l_1=\vert\Lambda_1\vert$ its cardinality and $\{A^{(1)}_i\}_{i\in\Lambda_1}$ the corresponding rational functions. This implies in particular that $\{A^{(1)}_i f_i^n\}_{i\in\Lambda_1}$ satisfy a minimal linear relation. As we have already shown earlier  (see~(\ref{N1Siegel}) and~(\ref{basisV1})), we have
$$
(r-1)N_1\leq nd + O(n/K)
$$
and there exists a basis $\v^{(1)}_1,\ldots,\v^{(1)}_{l_1-1}$ of $V_{\Lambda_1}$ satisfying
$$
h(\v^{(1)}_i)\leq N_1h(\Pt)+O\left(\frac{n}{K}h(\Pt)+nh(\Pt)^{1/2}+nK\right).
$$

We also set $J_1=\emptyset$ and $\varphi_1(i)=1$ for $i=1,\ldots,r$.\\

For $j\ge 2$ we go ahead similarly, but modifying somewhat the requirements for the $A_j$, taking into account the previous steps. More precisely, let $j\geq2$ and assume to have already constructed $N_1,\ldots,N_{j-1}$, $\Lambda_1,\ldots,\Lambda_{j-1}$ and $J_{j-1}$, $\varphi_{j-1}$. If $\{\Lambda_1,\ldots,\Lambda_{j-1}\}$ is connected and $\Lambda_1\cup\cdots\cup\Lambda_{j-1}=\{1,\ldots,r\}$ we put $s=j-1$ and we stop here the process.

Otherwise, we choose a subset $J_j$ of $\Lambda_1\cup\cdots\cup\Lambda_{j-1}$ such that $\vert J_j\cap C\vert=1$ for each 
connected component of $\{\Lambda_1,\ldots,\Lambda_{j-1}\}$. We need to show that we can choose $J_j$ disjoint from $\Lambda_1\cup\cdots \cup\Lambda_{j-2}\setminus J_{j-1}$. If the set $\Lambda_{j-1}$ does not intersect $\Lambda_1\cup\cdots\cup\Lambda_{j-2}$, the connected components of $\{\Lambda_1,\ldots,\Lambda_{j-1}\}$ are the connected components of $\{\Lambda_1,\ldots,\Lambda_{j-2}\}$ plus the set $\Lambda_{j-1}$ itself. Thus we may choose $J_j=J_{j-1}\cup\{i_0\}$ where $i_0$ is any element of $\Lambda_{j-1}$. 
If otherwise $\Lambda_{j-1}$ intersects $\Lambda_1\cup\cdots\cup\Lambda_{j-2}$, each connected component of $\{\Lambda_1,\ldots,\Lambda_{j-1}\}$ contains at least one connected component of $\{\Lambda_1,\ldots,\Lambda_{j-2}\}$. Thus we may choose $J_j$ as a suitable subset of $J_{j-1}$.

Then we let  $N_j$ be  the minimum of the set of integers $N\geq0$ such that there exist $\Lambda\subseteq\{1,\ldots,r\}$ and rational functions $A_i\in\F$ ($i\in\Lambda$) not all zero, satisfying\footnote{Note that we are prescribing somewhat more stringent conditions than before on the indices inside the subsets previously defined.}  
\begin{equation}
\label{sistema}
\begin{cases}
A_i\in L((N_{\varphi_{j-1}(i)}-1)Q),&\hskip-1cm \hbox{ if }Êi\in(\Lambda_1\cup\cdots\cup\Lambda_{j-1}\backslash J_j)\cap\Lambda,\\ 
A_i\in L(NQ),&\hskip-1cm \hbox{ if }Êi\in\Lambda \hbox{ and }Êi\not\in\Lambda_1\cup\cdots\cup\Lambda_{j-1}\backslash J_j,\\
h(A_i)\leq nK,&\hskip-1cm\hbox{ for } i\in\Lambda,\\
(A_if_i^n)_{i\in\Lambda} \hbox{ are linearly dependent.}
\end{cases}
\end{equation}

We remark that the set of such $N$ is indeed not empty, as we easily see since $\{1,\ldots,r\}\backslash(\Lambda_1\cup\cdots\cup\Lambda_{j-1}\backslash J_j)$ has cardinality $\geq2$ (for otherwise $\{\Lambda_1,\ldots,\Lambda_{j-1}\}$ would be  connected and $\Lambda_1\cup\cdots\cup\Lambda_{j-1}=\{1,\ldots,r\}$).

We select a minimal set $\Lambda_j$ among all sets $\Lambda$ which realize the minimum defining $N_j$. We denote by $\{A_i^{(j)}\}_{i\in\Lambda_j}$ the corresponding rational functions. Thus $\{A^{(j)}_i f_i^n\}_{i\in\Lambda_j}$ satisfy a minimal linear relation.

We finally set 
\begin{equation}
\label{varphi}
\varphi_j(i)=
\begin{cases}
\varphi_{j-1}(i) & \hbox{ if }Êi\in\Lambda_1\cup\cdots\cup\Lambda_{j-1}\backslash J_j,\\
j &  \hbox{otherwise, \ie  if }Êi\not\in\Lambda_1\cup\cdots\cup\Lambda_{j-1} \hbox{ or } i\in J_j.
\end{cases}
\end{equation}\

Thus, by the minimality of $N_j$, there are no subsets $\Lambda$ and rational functions $(A_i)_{i\in\Lambda}$ not all zero such that 
$$
\begin{cases}
A_i\in L((N_{\varphi_j(i)}-1)Q),&\hskip-1cm \hbox{ for} \ Êi\in\Lambda,\\ 
h(A_i)\leq nK,&\hskip-1cm\hbox{ for } i\in\Lambda,\\
(A_if_i^n)_{i\in\Lambda} \hbox{ are linearly dependent.}
\end{cases}
$$\

We notice that for $j=2,\ldots,s$ we have
\begin{equation}
\label{joker}
\varphi_{j-1}(i)=j-1\hbox{ for } i\in J_j,
\end{equation}
since $J_j$ is disjoint from $\Lambda_1\cup\cdots \cup\Lambda_{j-2}\setminus J_{j-1}$.\\

{\smallskip\noindent{\bf Proof of Claim~\ref{claim} i).}\quad}
Since at each step we have added conditions, we have 
\begin{equation}
\label{ordine}
N_1\leq N_2\leq\cdots \leq N_s
\end{equation}
as required.\\

{\smallskip\noindent{\bf Proof of Claim~\ref{claim} ii) and iii)}\quad} We first showÊthat for $j=2,\ldots,s$,  the set $\Lambda_j$ is contained in no connected component of $\{\Lambda_1,\ldots,\Lambda_{j-1}\}$. 

To verify this assertion, let us assume  by contradiction that $\Lambda_j$ is contained in a connected component $C$ of $\{\Lambda_1,\ldots,\Lambda_{j-1}\}$. Let $C\cap J_j=\{i_0\}$. Then for $i\in \Lambda_j$, $i\neq i_0$, we have $A_i^{(j)}\in L((N_{\varphi_{j-1}(i)}-1)Q)$. Using the minimal linear relation $\sum_{i\in\Lambda_j} a_iA^{(j)}_i f_i^n=0$ and the fact that $Q$ is not in the support of any of  the divisors $\div(f_i)$ we see that $A_{i_0}^{(j)}\in L(N^*Q)$ with $N^*=\max_i(N_{\varphi_{j-1}(i)}-1)$. By~(\ref{ordine}) we have 
$N^*\leq N_{j-1}$ and, by~(\ref{joker}), $j-1=\varphi_{j-1}(i_0)$. Thus $A_i^{(j)}\in L((N_{\varphi_{j-1}(i)}-1)Q)$ for {\sl all} $i\in\Lambda_j$. This contradicts the minimality in the definition of $N_{j-1}$ (see the remark after~({\ref{varphi})).

Thus for $j=2,\ldots,s$ the set $\Lambda_j$ is contained in no connected component of $\{\Lambda_1,\ldots,\Lambda_{j-1}\}$. This proves assertion ii) of Claim~\ref{claim} and ensures that the inductive process ends somewhere. Thus Claim~\ref{claim} iii) also holds by inductive construction.\\

We still have to check assertions iv) and v) of Claim~\ref{claim}. To prove assertion iv) we first need the following lemma.  Let, as in Claim~\ref{claim} iv), $t_1=\dim(V_{\Lambda_1})$ and
$$
t_j=\dim(V_{\Lambda_1}+\cdots+V_{\Lambda_j})-\dim(V_{\Lambda_1}+\cdots+V_{\Lambda_{j-1}}).
$$
for $j=2,\ldots,s$. 
\begin{lemma}
\label{combinatoriale}
~
\begin{itemize}
\item[i)]
For $j=2,\ldots,s$ we have:
$$
\vert\Lambda_1\cup\cdots\cup\Lambda_{j-1}\vert-\vert J_j\vert
=\dim(V_{\Lambda_1}+\cdots+V_{\Lambda_{j-1}})=\sum_{j'=1}^{j-1} t_{j'}
$$
\item[ii)]
Let $j$, $j'$ be two integers with $j\geq 2$ and $1\leq j'\leq j-1$. Then the set of $i\in\{1,\ldots,r\}$ such that $\varphi_j(i)=j'$ has cardinality $t_{j'}$.
\end{itemize}
\end{lemma}
\DIM The first assertion follows from Remark~\ref{increase} ii), since $J_j$ has cardinality equal to the number of connected components of $\Lambda_1\cup\cdots\cup\Lambda_{j-1}$.\\

We prove the second assertion by induction on $j$. For $j=2$ we have by construction $\#\{i\;\vert\; \varphi_2(i)=1\}=\vert \Lambda_1\vert-1=\dim(V_{\Lambda_1})=t_1$.

Let $j>2$ and assume that $\#\{i\;\vert\; \varphi_{j-1}(i)=j'\}=t_{j'}$ for $j'=1,\ldots,j-2$. We want to show that 
 $\#\{i\;\vert\; \varphi_j(i)=j'\}=t_{j'}$ for $j'=1,\ldots,j-1$. 
 
Assume first $j'\leq j-2$. Let $i$ be such that $\varphi_j(i)=j'$. Then $\varphi_j(i)\neq j$, thus by~(\ref{varphi}) $\varphi_j(i)=\varphi_{j-1}(i)$. This shows that $\{i\;\vert\; \varphi_j(i)=j'\}\subseteq\{i\;\vert\; \varphi_{j-1}(i)=j'\}$. On the other hand, let $i$ such that $\varphi_{j-1}(i)=j'$. Then $\varphi_{j-1}(i)\neq j-1$ and~(\ref{varphi}) (with $j$ replaced by $j-1$) shows that 
$$
i \in\Lambda_1\cup\cdots\cup\Lambda_{j-2}\backslash J_{j-1}\subseteq\Lambda_1\cup\cdots\cup\Lambda_{j-1}.
$$
Moreover $i\not\in J_j$ by~(\ref{joker}). Thus, by~(\ref{varphi}), $\varphi_j(i)=\varphi_{j-1}(i)$. This proves that 
$\{i\;\vert\; \varphi_{j-1}(i)=j'\}\subseteq\{i\;\vert\; \varphi_j(i)=j'\}$. Putting together the two inclusions we see that 
$$
\{i\;\vert\; \varphi_j(i)=j'\}=\{i\;\vert\; \varphi_{j-1}(i)=j'\}.
$$
Thus, by induction, $\#\{i\;\vert\; \varphi_j(i)=j'\}=\#\{i\;\vert\; \varphi_{j-1}(i)=j'\}=t_{j'}$.

Assume now $j'=j-1$. By~(\ref{varphi}) 
$$
\{i\;\vert\; \varphi_j(i)=j-1\}=\{i\in \Lambda_1\cup\cdots\cup\Lambda_{j-1}\backslash J_j\;\vert\;\varphi_{j-1}(i)=j-1\}
$$
By~(\ref{joker}) we have $\varphi_{j-1}(i)=j-1$ on $J_j$. 
By~(\ref{varphi}) (again with $j$ replaced by $j-1$) we still have $\varphi_{j-1}(i)=j-1$ outside $\Lambda_1\cup\cdots\cup\Lambda_{j-2}$ and thus, {\sl a fortiori} outside $\Lambda_1\cup\cdots\cup\Lambda_{j-1}$. These facts and the first assertion of the present lemma imply 
\begin{align*}
\#\{i\;\vert\; \varphi_j(i)=j-1\}
&=\#\{i\;\vert\;\varphi_{j-1}(i)=j-1\}-(r-\vert\Lambda_1\cup\cdots\cup\Lambda_{j-1}\vert)-\vert J_j\vert\\
&=\dim(V_{\Lambda_1}+\cdots+V_{\Lambda_{j-1}})-(r-\#\{i\;\vert\;\varphi_{j-1}(i)=j-1\}).
\end{align*}
By induction, 
\begin{align*}
r-\#\{i\;\vert\;\varphi_{j-1}(i)=j-1\}
&=\sum_{j'=1}^{j-2}\#\{i\;\vert\;\varphi_{j-1}(i)=j'\}
=\sum_{j'=1}^{j-2}t_{j'}\\
&=\dim(V_{\Lambda_1}+\cdots+V_{\Lambda_{j-2}}).
\end{align*}
Putting together the last two displayed equations, we get
$$
\#\{i\;\vert\; \varphi_j(i)=j-1\}=\dim(V_{\Lambda_1}+\cdots+V_{\Lambda_{j-1}})-\dim(V_{\Lambda_1}+\cdots+V_{\Lambda_{j-2}})=t_{j-1}
$$
as desired.
\CVD

{\smallskip\noindent{\bf Proof of Claim~\ref{claim} iv).}\quad
We proceed by induction. For $j=1$ assertion iv) was already proved in~(\ref{N1Siegel}). Let $j\geq 2$.} In a similar way as we have done for $N_1$, we are going to provide an upper bound for $N_j$ using Lemma~\ref{siegel-substitute}. Let $\e=c/K$ for a sufficiently large constant $c$ and choose $N$ as the smallest integer such that 
\begin{equation}
\label{Ndef2}
\sum_{j'=1}^{j-1}t_{j'}(N_{j'}-1)+\Big(r-1-\sum_{j'=1}^{j-1}t_{j'}-\e\Big)N\geq(1+\e)dn+c_0.
\end{equation} 
We notice that $N=O(n)$.
\begin{fact} We have 
\label{Ninequa}
$$
N_{j-1}\leq N.
$$
\end{fact}
\DIM By Claim~\ref{claim} iv) with $j$ replaced by $j-1$ (which holds by the present inductive assumption) we have 
$$
\sum_{j'=1}^{j-2}t_{j'}N_{j'}+\Big(r-1-\sum_{j'=1}^{j-2}t_{j'}\Big)N_{j-1}\leq nd +O(n/K).
$$ 
Thus 
\begin{align*}
\Big(r-1-&\sum_{j'=1}^{j-1}t_{j'}-\e\Big)(N_{j-1}-N)\\
&=\Big(r-1-\sum_{j'=1}^{j-2}t_{j'}\Big)N_{j-1}-t_{j-1}N_{j-1}-\e N_{j-1}-\Big(r-1-\sum_{j'=1}^{j-1}t_{j'}-\e\Big)N\\
&\leq \Big( nd +O(n/K)\Big)-\sum_{j'=1}^{j-2}t_{j'}N_{j'}-t_{j-1}N_{j-1}-\Big(r-1-\sum_{j'=1}^{j-1}t_{j'}-\e\Big)N\\
&\leq \Big( nd +O(n/K)\Big)-\Big((1+\e)dn+c_0)\Big)\\
&=\big(-c d +O(1)\big)\frac{n}{K}<0,
\end{align*}
if $c$ is a sufficiently large constant to kill the implicit constant in the $O(1)$.
\CVD
We are going to apply Lemma \ref{siegel-substitute} with $f_1,\ldots,f_r$ and with
$$
M_i=
\begin{cases}
N_{\varphi_{j-1}(i)}-1,& \hbox{ if }Êi\in(\Lambda_1\cup\cdots\cup\Lambda_{j-1}\backslash J_j);\\ 
N,& \hbox{ if }Êi\not\in\Lambda_1\cup\cdots\cup\Lambda_{j-1}\backslash J_j
\end{cases}
$$
($i=1,\ldots,r$). Since $Q$ is neither a zero nor a pole of $f_1,\ldots,f_r$ and since 
$N_1\leq \cdots \leq N_{j-1}\leq N$ (by Claim~\ref{claim} i) and by Remark~\ref{Ninequa}), we have 
\begin{equation}
\label{E}
M=\max_iM_i=N.
\end{equation}
We recall that, by Lemma~\ref{combinatoriale} i),
$$
\vert \Lambda_1\cup\cdots\cup\Lambda_{j-1}\vert-\vert J_j\vert =\sum_{j'=1}^{j-1} t_{j'}
$$
and, by~(\ref{varphi}) and by Lemma~\ref{combinatoriale} ii),
$$
\#\{i\in(\Lambda_1\cup\cdots\cup\Lambda_{j-1}\backslash J_j)\;\vert\; \varphi_{j-1}(i)=j'\}=
\#\{i\;\vert\; \varphi_j(i)=j'\}=t_{j'}
$$
for $j'=1,\ldots,j-1$. Thus 
\begin{equation}
\label{U}
S=\sum_{i=1}^rM_i=\sum_{j'=1}^{j-1}t_{j'}(N_{j'}-1)+\Big(r-\sum_{j'=1}^{j-1}t_{j'}\Big)N.
\end{equation} 
By~(\ref{Ndef2}),~(\ref{E}) and~(\ref{U}), $S\geq (1+\e)(N+dn)+c_0$ and the Dirichlet exponent $\varrho$ of Lemma~\ref{siegel-substitute} satisfies
$$
\varrho=\frac{N+dn}{S-N-dn-c_0}\leq\frac{1}{\e}.
$$
By that lemma (if $c$ is a sufficiently large constant) there exist rational functions $A_1,\ldots, A_r$ not all zero and of height $\leq nK$ such that 
$$
\begin{cases}
A_i\in L((N_{\varphi_{j-1}(i)}-1)Q),& \hbox{ if }Êi\in(\Lambda_1\cup\cdots\cup\Lambda_{j-1}\backslash J_j);\\ 
A_i\in L(NQ),& \hbox{ if }Êi\not\in\Lambda_1\cup\cdots\cup\Lambda_{j-1}\backslash J_j
\end{cases}
$$
satisfying 
$$
A_1f_1^n+\ldots+A_rf_r^n=0.
$$
Not all $A_1,\ldots,A_r$ are zero,  and we see that the non-zero ones among $A_1f_1^n,\ldots,A_rf_r^n$ sum up to zero, and so they are linearly dependent. Choosing $\Lambda$ as the set of $i$ such that $A_i\neq0$ we see that~(\ref{sistema}) is satisfied. By minimality of $N_j$ we have $N_j\leq N$. Since $N$ is the smallest integer satisfying~(\ref{Ndef2}) and since $\e=c/K$, we have 
$$
\sum_{j'=1}^{j-1}t_{j'}N_{j'}+\Big(r-1-\sum_{j'=1}^{j-1}t_{j'}\Big)N_j\leq nd+O(n/K).
$$
as required.
This conclude the proof of assertion iv) of Claim~\ref{claim}.\\

{\smallskip\noindent{\bf Proof of Claim~\ref{claim} v).}\quad
For $j=1$ assertion v) was already proved in~(\ref{basisV1}). Let $j\geq 2$. }As we did 
for $j=1$, we apply Lemma~\ref{Wronskiano} to find a suitable basis of $V_{\Lambda_j}$. Let, for $i\in\Lambda_j$, 
$$
M'_i=
\begin{cases}
N_{\varphi_{j-1}(i)}-1,&\hbox{ if }Êi\in(\Lambda_1\cup\cdots\cup\Lambda_{j-1}\backslash J_j)\cap\Lambda_j,\\ 
N_j,&\hbox{ if }Êi\in\Lambda_j \hbox{ and }Êi\not\in\Lambda_1\cup\cdots\cup\Lambda_{j-1}\backslash J_j
\end{cases}
$$
and $M'=\max_{i\in\Lambda_j} M'_i=N_j$. Let also, as in Lemma~\ref{Wronskiano}, 
$$
\Theta=\max\Big(1,\sum_{i\in\Lambda_j}M'_i-(M'+nd_j)\Big)=\max\Big(1,\sum_{i\in\Lambda_j}M'_i-(N_j+nd_j)\Big)
$$ 
with $d_j=d_{\Lambda_j}=-\deg\div(f_i)_{i\in\Lambda_j}$. We contend that  $\Theta=O(n/K)$. As for $j=1$, we prove this assertion using Lemma~\ref{siegel-substitute}, now for the $f_i~(i\in\Lambda_j)$ with $M_i$ the same as $M'_i$ except for $N_j-1$ in place of $N_j$, \ie
$$
M_i=
\begin{cases}
N_{\varphi_{j-1}(i)}-1,&\hbox{ if }Êi\in(\Lambda_1\cup\cdots\cup\Lambda_{j-1}\backslash J_j)\cap\Lambda_j,\\ 
N_j-1,&\hbox{ if }Êi\in\Lambda_j \hbox{ and }Êi\not\in\Lambda_1\cup\cdots\cup\Lambda_{j-1}\backslash J_j.
\end{cases}
$$
Let $S=\sum_{i\in\Lambda_j} M_i=\sum_{i\in\Lambda_j} M'_i+O(1)$ and $M=\max_{i\in\Lambda_j} M_i=N_j-1=O(n)$.  The Dirichlet exponent $\varrho$ is then
$$
\varrho=\frac{N_j-1+nd_j}{S -N_j-nd_j-c_0+1}.
$$
By Lemma~\ref{siegel-substitute}, there exist rational functions $B_i$ not all zero  such that $B_i\in L(M_iQ)$ for $i\in\Lambda_j$,  $\sum_{i\in\Lambda_j} B_if_i^n=0$ and 
$$
h(B_i)= O((\varrho+1) n).
$$
By the minimality of $N_j$, we cannot have $\max_i h(B_i)\leq nK$. Thus $\varrho\geq K/c$, where $c$ is a sufficiently large constant to kill the implicit constant in the last $O()$. This implies
$$
\sum_{i\in\Lambda_j}M'_i -(N_j+nd_j)\leq \frac{c}{K}(N_j+nd_j)+O(1)=O(n/K)
$$
as required.\\

We apply Lemma~\ref{Wronskiano} to the rational functions $\{A^{(j)}_i\}_{i\in\Lambda_j}$, taking into account:
$$
h(\Pt) \geq 1,\quad h(A_i)\leq nK,\quad M'_i=N_j=O(n),\quad \Theta=O(n/K).
$$
By this lemma, there exists a basis $\v^{(j)}_1,\ldots,\v^{(j)}_{l_j-1}$ of $V_{\Lambda_j}$ satisfying
\begin{align*}
h(\v^{(j)}_i)
&\leq M' h(\Pt)+O\left(\Theta (h(\Pt)+1)+(n+M')(1+h(\Pt)^{1/2})+\max h(A_i)\right)\\
&= N_jh(\Pt)+O\left(\frac{n}{K}h(\Pt)+nh(\Pt)^{1/2}+nK\right).
\end{align*}
This proves assertion v) of Claim~\ref{claim}. 

\section{Proof of Theorem~\ref{specialization}}
\label{coro}

Let $\Gamma\subset\Gm^r(\F)$ be a finitely generated constant-free  subgroup, and let $V$ be an algebraic subvariety of $\Gm^r$, defined over $\F$. By writing $V$ as an intersection of hypersurfaces, we see that it is enough to prove Theorem~\ref{specialization}  and its addendum for a hypersurface $V$. We may further assume that $V$ is a hyperplane. The case of a general hypersurface can indeed be easily deduced using an isogeny. For technical reasons, it is convenient to homogenize the statement in the linear case:
\begin{proposition}
\label{linear}
Let $\Gamma\subset\Gm^r(\F)$ be a finitely generated subgroup such that 
\begin{equation}
\label{H}
\forall\fb\in\Gamma,\quad \forall i,j=1,\ldots,r,\quad f_i/f_j\in\Qb \Longrightarrow f_i/f_j\in\Qb^*_{\rm tors}.
\end{equation}
Let $\gammag=(\gamma_1,\ldots ,\gamma_r)\in\Gamma$ and $\theta_1,\ldots,\theta_r\in\F$ such that $\theta_1\gamma_1+\cdots + \theta_r\gamma_r\neq0$. Then the height of $\Pt\in\Cu({\overline{\mathbb Q}})$ such that 
\begin{equation}
\label{hyperplane}
\theta_1(\Pt) \gamma_1(\Pt)+\cdots + \theta_r(\Pt)\gamma_r(\Pt)=0
\end{equation}
is bounded from above in terms only on $\Gamma$ and $\theta_1,\ldots,\theta_r$, independently of $\gamma_1,\ldots,\gamma_r$. The same conclusion holds without the assumption~\eqref{H}, if $\theta_1,\ldots,\theta_r\in\Qb$ and $\Gamma/\Gamma\cap\Gm^r(\Qb)$ is of rank $1$.
\end{proposition}
\DIM Replacing $\{1,\ldots,r\}$ by a subset (and $\Gamma$ by its projection on the coordinates in the subset) we may assume that there are no proper vanishing subsums in~\eqref{hyperplane}. Dividing~\eqref{hyperplane} by $\gamma_1(\Pt)$ (and replacing $\Gamma$ by $\{(1,f_2/f_1,\ldots,f_r/f_1)\;\vert\; \fb\in\Gamma\}$) we may also assume $\gamma_1=1$ and 
\begin{equation}
\label{H2}
\Gamma\subseteq \{x_1=1\}\;.
\end{equation}\

Assume first that $\Gamma/\Gamma\cap\Gm^r(\Qb)$ is of rank $1$ and $\theta_i\in\Qb$. Thus there exists $\fb=(1,f_2,\ldots,f_r)\in\Gamma$ such that $\Gamma=\Gamma\cap\Gm^r(\Qb)\oplus\langle\fb\rangle$. Then $\gammag=(1,c_2f_2^n,\ldots,c_rf_r^n)$ for some $c_i\in\Qb$. We may assume that there exists some $i>1$ such that $f_i$ is non constant, since otherwise $\theta_1\gamma_1+\cdots + \theta_r\gamma_r$ is a non zero constant and equation~\eqref{hyperplane} does not have solutions. Thus we can apply Theorem~\ref{main}, and we find that the solutions of~\eqref{hyperplane} have bounded height.\\

Assume now that $\Gamma$ satisfies~\eqref{H}. Then $\Gamma=\Gamma_{\rm tors}\oplus\Gamma'$ where $\Gamma'$ is freely generated by, say, $\gb_1,\ldots,\gb_\kappa$. We now use Dirichlet's Theorem in a way which is inspired by a method appearing already in Bombieri's paper~\cite{Bo}, especially Lemma 4 therein. There exist $\omegag\in\Gamma_{\rm tors}$ and integers $\lambda_1,\ldots,\lambda_\kappa$ such that 
$$
\gammag=\omegag\gb_1^{\lambda_1}\cdots\gb_\kappa^{\lambda_\kappa}.
$$
Let $A=\max\vert\lambda_j\vert$. Since the equation~(\ref{hyperplane}) is not trivial (and since $\Gamma_{\rm tors}$ is finite), we can obviously assume $A$ unbounded. Let $Q\geq 1$ be an integer which will be fixed later, independently of $A$. By Dirichlet's Theorem on simultaneous approximation, there exists a positive integer $q\leq Q^{\kappa}$ and integers $p_j$ such that 
$$
\left\vert q\frac{\lambda_j}{A}-p_j\right\vert< \frac{1}{Q}.
$$
Let $n$ be the integral part of $A/q$. We write $\lambda_j=np_j+r_j$ and we set, for $i=1,\ldots,r$, 
$$
\rhog= \prod_{j=1}^\kappa \gb_j^{r_j}\in\Gamma,
\qquad \fb=Ê\prod_{j=1}^\kappa \gb_j^{p_j}\in\Gamma,
\qquad \alpha_i=\omega_i\theta_i(P)\rho_i(P).
$$
Since $\gamma_i=\omega_i\rho_if_i^n$, equation~\eqref{hyperplane} can be rewritten as $\alpha_1 f_1(\Pt)^n+\cdots+\alpha_r f_r(\Pt)^n=0$ (without proper vanishing subsums).

We have
$$
\vert p_j\vert\leq \left\vert q\frac{\lambda_j}{A}\right\vert+ \frac{1}{Q}\leq q+Q^{-1}\leq 2Q^{\kappa}.
$$
Thus, $\fb$ belongs to a finite set, depending only on the generators $\gb_1,\ldots,\gb_\kappa$ and on $Q$, as the exponents $\lambda_j$ vary. Moreover
$$
\vert r_j\vert=\left\vert\frac{A}{q}\left(q\frac{\lambda_j}{A}-p_j\right)-\left(n-\frac{A}{q}\right)p_j\right\vert
\leq (n+1)Q^{-1}+2Q^{\kappa}.
$$\\
Thus $d(\theta_i\rho_i)= O(n/Q+Q^{r\kappa})$, where the implicit constant in the big-$O$ depends only on $\theta_1,\ldots,\theta_r$ and on $\gb_1,\ldots,\gb_\kappa$.\\
We recall that $\alpha_i=\omega_i\theta_i(P)\rho_i(P)$, with $\omega_i$ root of unity. Using Lemma \ref{HeightMachine} (for the two functions $1$ and $\theta_i\rho_i$), we see easily that
$$
h(\alphag)\leq h(\alpha_1)+\cdots+h(\alpha_r)\leq C_1(n/Q+Q^{\kappa})h(\Pt)
$$
where the constant $C_1$ depends only on $\theta_1,\ldots,\theta_r$ and on $\gb_1,\ldots,\gb_\kappa$ (and neither on $Q$ nor on $A$). We now choose $Q=[4rC_1]+1$.

There is a $j$ such that $\lambda_j=\pm A$ and thus some of the $p_j$ are not zero. Since $g_1,\ldots,g_\kappa$ is a basis of $\Gamma'$ we have $\fb\not\in\Gamma_{\rm tors}$.  By~\eqref{H2} we have $f_1=1$. Thus, by~\eqref{H}, 
$f_i=f_i/f_1$ is non-constant for some $i>1$ and we can apply Theorem~\ref{main}. 

Let $C$ be the constant appearing in this theorem, which depends on $f_1,\ldots,f_r$ and thus only on $\theta_1,\ldots,\theta_r$ and on $\gb_1,\ldots,\gb_\kappa$  (since the rational functions $f_i$ belong to a finite set, depending only on the generators $\gb_1,\ldots,\gb_\kappa$ and on $Q$ and $Q$ has already been fixed as $Q=[4rC_1]+1$). 
Since $A$ is unbounded and $Q$ does not depend on $A$, we may assume $n=[A/q]\geq C$ and $n\geq 4rC_1Q^\kappa$. By Theorem~\ref{main} and by the inequalities $h(\alphag)\leq C_1(n/Q+Q^{\kappa})h(\Pt)$, $Q\geq 4rC_1$ and $n\geq 4rC_1Q^\kappa$, 
$$
h(\Pt)\leq \frac{r h(\alphag)}{n}+C \leq rC_1(1/Q+Q^{\kappa}/n)h(\Pt)+C\leq h(\Pt)/2+C\;.
$$
We deduce that $h(\Pt)$ is bounded, as claimed.
\CVD

\end{document}